\documentclass[12pt]{article}
\usepackage{amsthm, amsmath, amsfonts}
\usepackage{geometry}
\usepackage{graphicx}
\usepackage{color}
\usepackage{hyperref}
\usepackage{comment}
\usepackage{bbm}

\newtheorem{theorem}{Theorem}

\newtheorem{lemma}[theorem]{Lemma}
\newtheorem{proposition}[theorem]{Proposition}

\newtheorem{remark}[theorem]{Remark}

\theoremstyle{remark}

\theoremstyle{definition}
\newtheorem{definition}[theorem]{Definition}

\def\cS{\mathcal{S}}

\def\cQ{\mathcal{Q}}

\def\cI{\mathcal{I}}

\def\cG{\mathcal{G}}
\def\cF{\mathcal{F}}

\def\cS{\mathcal{S}}

\newcommand{\C}{\mathbbm{C}}
\newcommand{\D}{\mathbbm{D}}
\newcommand{\E}{\mathbbm{E}}
\newcommand{\N}{\mathbbm{N}}

\newcommand{\R}{\mathbbm{R}}
\renewcommand{\P}{\mathbbm{P}}

\newcommand{\wt}{\widetilde}
\newcommand{\wh}{\widehat}
\newcommand{\sgn}{\operatorname{sign}}
\newcommand{\rta}{\rightarrow}

\newcommand{\ceil}[1]{\left\lceil #1 \right\rceil}

\newcommand{\1}{\mathbf{1}}
\newcommand{\eps}{\varepsilon}
\def\e{\mathbf{e}}
\def\a{\mathbf{a}}
\def\at{\wt{\mathbf{a}}}
\def\x{\mathbf{x}}
\def\xt{\wt{\mathbf{x}}}
\def\w{\mathbf{w}}
\def\wtt{\wt{\mathbf{w}}}
\def\bad{\Xi_{\rm bad}}

\newcommand{\polylog}{\operatorname{polylog}(n)}
\newcommand{\Linfm}[1]{\left\| #1 \right\|_{m,\infty}}

\DeclareMathOperator*{\argmin}{arg\,min}

\def\@rst #1 #2other{#1}
\newcommand\MR[1]{\relax\ifhmode\unskip\spacefactor3000 \space\fi
  \MRhref{\expandafter\@rst #1 other}{#1}}
\newcommand{\MRhref}[2]{\href{http://www.ams.org/mathscinet-getitem?mr=#1}{MR#2}}

\def\MR#1{\href{http://www.ams.org/mathscinet-getitem?mr=#1}{MR#1}}

\begin{document}

\author{
\begin{tabular}{c}Nina Holden\thanks{ETH Z\"urich, nina.holden$@$eth-its.ethz.ch}\end{tabular}\;
\begin{tabular}{c}Robin Pemantle\thanks{University of Pennsylvania, pemantle$@$math.upenn.edu}\end{tabular}\;
\begin{tabular}{c}Yuval Peres\thanks{yperes$@$gmail.com}\end{tabular}
\begin{tabular}{c}Alex Zhai \end{tabular}}
\setcounter{tocdepth}{2}
\title{
	Subpolynomial trace reconstruction for random strings and arbitrary deletion probability
}
\date{}
\maketitle



\begin{abstract}
  The insertion-deletion channel takes as input a bit string $\x\in
  \{0,1\}^{n}$, and outputs a string where bits have been deleted and
  inserted independently at random. The trace reconstruction problem
  is to recover $\x$ from many independent outputs (called ``traces'')
  of the insertion-deletion channel applied to $\x$. We show that if
  $\x$ is chosen uniformly at random, then $\exp(O(\log^{1/3} n))$
  traces suffice to reconstruct $\x$ with high probability. For the
  deletion channel with deletion probability $q < 1/2$ the earlier upper
  bound was $\exp(O(\log^{1/2} n))$. The case of $q \geq 1/2$ or the
  case where insertions are allowed has not been previously analyzed,
  and therefore the earlier upper bound was as for worst-case strings,
  i.e., $\exp(O( n^{1/3}))$. We also show that our reconstruction algorithm runs in $n^{1+o(1)}$ time.

  A key ingredient in our proof is a delicate two-step alignment
  procedure where we estimate the location in each trace corresponding
  to a given bit of $\x$. The alignment is done by viewing the strings
  as random walks and comparing the increments in the walk associated
  with the input string and the trace, respectively.
\end{abstract}

\section{Introduction}

\begin{figure}[ht!]
  \centering
  \includegraphics[scale=0.95]{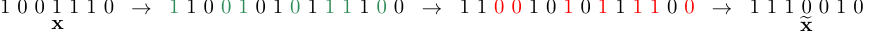}
  \caption{We obtain a trace $\xt$ by sending $\x$ through the
    insertion-deletion channel. Inserted bits are shown in green and
    deleted bits are shown in red.}
  \label{fig:trace}
\end{figure}

Learning a parameter from a sequence of noisy observations is a basic
problem in statistical inference and machine learning. The amount of
data required (known as the {\em sample complexity}) to learn the
parameter is of fundamental interest. A natural problem in this class
where the missing parameter is a bit string and it is unknown whether
the sample complexity is polynomial, is the trace reconstruction
problem for the {\bf insertion-deletion channel}. This channel takes
as input a string $\x=(x_0,x_1,\dots,x_{n-1}) \in \{0,1\}^{n}$ and
outputs a noisy version of it, where bits have been randomly inserted
and deleted. Let $q\in[0,1)$ be the deletion probability and let
  $q'\in[0,1)$ be the insertion probability. First, for each $j$,
    before the $j$th bit of $\x$ we insert $G_j-1$ uniform and
    independent bits, where the independent geometric random variables
    $G_j \ge 1 $ have parameter $1-q'$. Then we delete each bit of the
    resulting string independently with probability $q$.  The output
    string $\xt$ is called a {\bf trace}.  An example is shown in
    Figure~\ref{fig:trace}.

Suppose that the input string $\x$ is unknown. The {\bf trace
  reconstruction problem} asks the following: How many i.i.d.\ copies
of the trace $\xt$ do we need in order to determine $\x$ with high
probability? (A more formal problem description will be given in
Section \ref{sec:reconstruction}.)

There are two variants of this problem: the ``worst case'' and the
``average case'' (also referred to as the ``random case''). In the
worst case variant, we want to obtain bounds which hold uniformly over
all possible input strings $\x$. In the average case variant, the
input string is chosen uniformly at random.

In this paper, we study the average case. Holenstein, Mitzenmacher,
Panigrahy, and Wieder \cite{HMPW08} gave an algorithm for
reconstructing random strings from the deletion channel using
polynomially many traces, assuming the deletion probability $q$ is
sufficiently small. Peres and Zhai \cite{PZ17} proved that $\exp(O(
\log^{1/2}n ))$ many traces suffice for the deletion channel when the
deletion probability $q$ is below $1/2$. For $q \geq 1/2$, the
previous best bound was the same as for worst case strings, i.e.,
$\exp(O( n^{1/3} ))$, see \cite{DOS16,NP16}.

The following theorem is our main result, which improves the upper
bound for all $q \in [0,1)$ and also holds when we allow insertions. In particular, we answer the part of the first open question in \cite[Section 9]{M09} which concerns random strings.
\begin{theorem}\label{thm:main}
  For $n\in\N$ let $\x\in\{0,1 \}^n$ be a bit string where the bits
  are chosen uniformly and independently at random.  Given
  $q,q'\in[0,1)$ there exists $M>0$ such that for all $n$ we can
    reconstruct $\x$ with probability $1-o_n(1)$ using
    $\lceil\exp(M\log^{1/3}n)\rceil$ traces from the
    insertion-deletion channel with parameters $q$ and $q'$. Moreover,
    our algorithm runs in $n^{1 + o(1)}$ time.
\end{theorem}

An earlier version of this work appeared in COLT \cite{HPP18}; in this
updated version, we have simplified the proof and added an analysis of
the algorithm's running time.

We remark that the trace reconstruction
problem is significantly more difficult for $q>1/2$
  \footnote{
    Suppose that $q>1/2$, the string $\w$ is an arbitrary string of
    length $(1-q)n$, and $\x$ is a random string of length $n$. Then
    it holds with probability at least $1-\exp(-cn)$ that $\w$ is a subsequence
    of $\x$. To see this, observe that the number of bits in $\x$
    until we see $w_0$ is a geometric random variable of mean
    2. Iterating, existence (with probability $1-\exp(-cn)$) of an
    appropriate subsequence holds by concentration for the sum of
    independent geometric random variables of mean 2. Therefore, by a
    union bound, if $q>1/2$, then any subexponential collection of
    strings of length $(1-q)n$ (typical for traces of the deletion
    channel) are, with high probability, all substrings of a random
    string $\x$ of length $n$.
  }
and that the alignment algorithm used by Peres and Zhai fails
fundamentally in this case. Moreover, the upper bound
$\exp(O(\log^{1/3} n))$ in Theorem \ref{thm:main} is the best one can
obtain without also improving the upper bound $\exp(O(n^{1/3}))$ for
worst case strings. Indeed, given an arbitrary string of length $m =
\log_{2+\eps} n$ for $\eps>0$, this string will appear in a random
length $n$ string with probability converging to 1 as $n\rta\infty$.
In particular, a given worst case string of length $m$ is likely to
appear in our random string, and the best known algorithm for
reconstructing this string requires $\exp(\Omega(m^{1/3})) =
\exp(\Omega(\log^{1/3} n))$ traces. See Lemma 10 in \cite{MPV14} for
the details of this reduction.

We note also that our methods can be adapted easily to certain other
reconstruction problems, e.g., to the case where one allows
substitutions in addition to deletions and insertions, and the case
where the bits in the input $\x$ are independent Bernoulli($r$) random
variables for arbitrary $r\in(0,1)$, instead of $r=1/2$. There is also
a simple reduction (described in e.g.\ \cite{MPV14} and \cite{DOS16})
of the trace reconstruction problem for larger alphabets to the case
of bits. Moreover, as shown in \cite{MPV14}, trace reconstruction
becomes much easier if the alphabet size grows as $\Omega(\log n)$.

In Section \ref{sec:literature} we present some background and literature on the trace reconstruction problem, before we give a precise definition of the trace reconstruction problem in Section \ref{sec:reconstruction}. We give an outline of the proof of Theorem \ref{thm:main} in Section \ref{sec:outline} and we present some notation in Section \ref{sec:notation}. In Sections \ref{sec:test} to \ref{sec:bit-test} we prove the various ingredients which are needed for the proof of Theorem \ref{thm:main}, and in Section \ref{sec:proof-conclusion} we conclude the proof.

\section{Related work}
\label{sec:literature}

The trace reconstruction problem dates back to the early 2000's \cite{L01a,L01b,BKKM04}. Batu, Kannan, Khanna, and McGregor,
who were partially motivated by the study of genetic mutations,
considered the case where the deletion probability $q$ is decreasing
in $n$. They proved that if the original string $\x$ is random and the
deletion probability $q=O(1/\log n)$, then $\x$ can be constructed
with high probability using $O(\log n)$ samples. Furthermore, they
proved that if $q=O(n^{-(1/2+\eps)})$, then every string $\x$ can be
reconstructed with high probability with $O(n\log n)$ samples.

Holenstein, Mitzenmacher, Panigrahy, and Wieder \cite{HMPW08}
considered the case of random strings and constant deletion
probability. They gave an algorithm for reconstruction with
polynomially many traces when the deletion probability $q$ is less
than some small threshold $c$. The threshold $c$ is not given
explicitly in the work of \cite{HMPW08}, but was estimated in
\cite{PZ17} to be at most $0.07$.

The result of \cite{HMPW08} was improved by \cite{PZ17}. They showed
that a subpolynomial number of traces $\exp(O(\log^{1/2} n))$ is
sufficient for reconstruction, and they extended the range of allowed
$q$ to the interval $[0,1/2)$.

Our work improves the above results in three ways. First, we improve
the upper bound to $\exp(O(\log^{1/3} n))$. Second, we allow for any
deletion and insertion probabilities in $[0,1)$. Third, unlike
  \cite{PZ17}, our method works not only for the deletion channel, but
  also for the case where we allow insertions and substitutions.

It is shown by \cite{HMPW08} that $\exp(O(n^{1/2} \log n))$ traces
suffice for reconstruction with high probability with worst case
input. This was improved to $\exp(O(n^{1/3}))$ independently by De, O'Donnell, and Servedio
\cite{DOS16} and by Nazarov and Peres \cite{NP16}. Until the current work, the average
case upper bound was equal to the worst case upper bound for $q\geq
1/2$. The techniques developed by \cite{DOS16,NP16} are applied in the
current work and the work of \cite{PZ17} to certain shorter substrings
of our random string.

A lower bound of $\Omega(\log^2 n)$ was obtained in the average case in \cite{MPV14}, and a lower bound of $\Omega(n)$ was obtained in the worst case in \cite{BKKM04}. These bounds were improved to $\Omega(\log^{9/4}n/\sqrt{\log\log n})$ and $\Omega(n^{5/4}/\sqrt{\log n})$, respectively, by Holden and Lyons \cite{HL20}, and further to $\Omega(\log^{5/2} n/(\log\log n)^7)$ and $\Omega(n^{3/2}/\log^7 n)$ by Chase \cite{chase}.
Trace reconstruction for the
setting which allows insertions and substitution in addition to
deletions was considered in \cite{KM05,VS08,DOS16,NP16}. We refer to the introduction of 
\cite{DOS16} and the survey \cite{M09} for further background on
the deletion channel.

\section{The trace reconstruction problem} \label{sec:reconstruction}

To simplify notation, throughout the paper we will implicitly pad any
finite-length bit strings with infinitely many zeroes on the
right. Thus, expressions such as $\x(i:j)$ are well-defined for $i<j$
even if $j$ is larger than the length of $\x$. Let $\N = \{0,1,\dots
\}$, and let $\cS := \{ 0 , 1 \}^\N$ denote the space of infinite
sequences of zeroes and ones. We denote elements of $\cS$ by $\x :=
(x_0 , x_1 , \ldots)$. If $I\subset\R$ we will sometimes write $\x(I)$ instead of $\x(I\cap\N)$, and we use a similar convention for functions which are defined on (subsets of) $\N$.

Fix a deletion probability $q$ and an insertion probability $q'$ in
$[0,1)$, and let $p=1-q$ and $p'=1-q'$. We construct $\xt$ from $\x$
  by the procedure described above, i.e., first, for each $j\in\N$ we
  insert $G_j-1$ uniform and independent bits before the $j$th bit of
  $\x$. The geometric random variables $G_j$ are independent and
  satisfy
\[ \P[ G_j=v ]=(q')^{v-1}(1- q'),\qquad\forall v\in\{1,2,\dots \}. \]
Then we delete each bit of the resulting string independently with
probability $q$.

Let $\mu$ be the law of i.i.d.\ Bernoulli random variables with
parameter $1/2$.  We write $\P_\x := \P_{\delta_\x}$ for the law of
$\xt$ when $\x$ is fixed and write $\P := \P_\mu$ for the law of $\xt$
when $\x$ is picked according to $\mu$. We call the string $\xt$ a
trace. An example is given in Figure~\ref{fig:trace}.

\subsubsection*{Worst case reconstruction problem}

Let $q,q'\in[0,1)$. For any $N\in\N$ let $\P_{\x}^{N}$ denote the
probability measure associated with $N$ independent outputs of the
insertion-deletion channel $\P_{\x}$ with deletion (resp.\ insertion)
probability $q$ (resp.\ $q'$). For $n\in\N$ and $\x\in\{0,1 \}^n$ let
$\mathfrak{X}$ denote a collection of $N_n\in\N$ traces sampled
independently at random. We say that worst case strings of length $n$
can be reconstructed with probability $1 - o_n(1)$ from $N_n$ traces,
if there is a function $G\colon \cS ^{N_n} \to \{0,1 \}^n$, such that
for all $\x\in\cS$,
\[ \P_{\x}^{N_n}[G(\mathfrak X ) = \x(0:n-1)] = 1-o_n(1). \]

\subsubsection*{Average case reconstruction problem}

Let $\mu_n$ denote uniform measure on $\{0,1 \}^n$.  We say that
uniformly random strings of length $n$ can be reconstructed with
probability $1-o_n(1)$ from $N_n$ traces if we can find a set
$\cS_n\subset\{0,1 \}^n$ with $\mu_n(\cS_n)=1-o_n(1)$, and a function
$G\colon \cS^{N_n} \to \{0,1 \}^n$, such that for all $\x\in\cS$ for
which $\x(0:n-1) \in\cS_n$, we have
\[ \P_{\x}^{N_n}[G(\mathfrak X ) = \x(0:n-1)] = 1-o_n(1). \]
In particular, Theorem \ref{thm:main} says that uniformly random
strings can be reconstructed from $N_n :=
\lceil\exp(M\log^{1/3}n)\rceil$ traces with probability $1-o_n(1)$.

\section{Outline of proof}
\label{sec:outline}
We give here an informal description of the algorithm used to achieve
the bound in Theorem \ref{thm:main}. The bits of $\x$ will be recovered
one by one: for any $k,n \in \N$ with $k < n$ we assume $\x(0:k)$ is
already known, and we show that with probability $1-O(n^{-2})$ we can
use $\lceil\exp(M\log^{1/3}n)\rceil$ traces to determine the
subsequent bit $x_{k+1}$. We will reuse these same traces for each
step (i.e., for all values of $k$). Even with this reuse, by a
union bound, the reconstruction will succeed at every step with high
probability.

Three ingredients are required, as follows:
\begin{enumerate}
\item A Boolean test $T(\w, \wtt)$ on pairs $(\w, \wtt)$ of bit
  strings indicating whether $\wtt$ is a plausible match for the
  string $\w$ sent through the insertion-deletion channel.

\item An alignment procedure that uses the test $T$ repeatedly to
  produce for each of the independent traces $\xt$ an estimate $\tau$
  for the position in $\xt$ corresponding to the $k$-th bit of $\x$.

\item A bit recovery procedure based on a method
  of~\cite{PZ17,DOS16,NP16} to produce from the approximately aligned
  traces an estimate of the subsequent bit or bits.
\end{enumerate}

The argument of \cite{PZ17} follows the same overall structure, with
an alignment step followed by a reconstruction step for each bit in
the original string. However, the greedy alignment step of \cite{PZ17}
relies crucially on the assumption that the deletion probability $q <
1/2$, and that no insertions are allowed. We overcome this problem by
introducing a new kind of test for the alignment based on studying
correlations between blocks in the input string and in the trace.

\subsection{The test $T$}
We describe here a simplified version of the test $T(\w, \wtt)$, which
returns $1$ if there is a likely match and $0$ otherwise. Assume for
simplicity that $\w$ and $\wtt$ have the same length and that $q = q'$,
so the expected output length of the insertion-deletion channel is the
same as the input length. The test involves two parameters:
the length $\ell$ of the strings to test and another parameter
$\lambda \le \sqrt{\ell}$. We subdivide both $\w$ and $\wtt$ into
roughly $\ell/\lambda$ blocks of size $\lambda$. We use
$T_{\ell,\lambda}$ to denote the test using these parameters.

\begin{figure}
	\centering
	\includegraphics[scale=1]{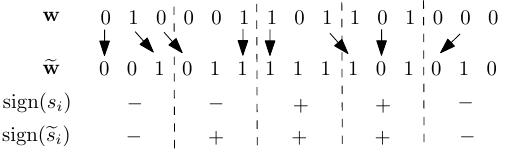}
	\caption{Illustration of the test $T_{\ell,\lambda}(\w,\wt{\w})$ with $\ell=15$ and $\lambda=3$. The arrows indicate that a bit in $\w$ was copied to a bit in $\wt\w$.}
	\label{fig:test}
\end{figure}

For each $i$, let $s_i$ be the number of $1$'s minus the number of
$0$'s in the $i$-th block of $\w$, and similarly define $\wt s_i$ for
$\wtt$. For some fixed constant $c > 0$ to be specified later, we
declare that
\[ T_{\ell,\lambda}(\w, \wtt) = \begin{cases}
  1 & \text{if } \displaystyle \sum_{i=1}^{\ceil{\ell/\lambda}} \sgn(s_i)
  \cdot \sgn(\wt s_i) > c \cdot \frac{\ell}{\lambda},\\ 0 &
  \text{otherwise}.
\end{cases} \]
See Figure \ref{fig:test} for an illustration.
The idea here is that if the bits in the $i$th block of $\wtt$ did not come from the $i$th block of $\w$,
then they will be independent random bits, and so each $\sgn(s_i)
\cdot \sgn(\wt s_i)$ will be $1$ or $-1$ with equal probability (if we ignore the small probability event on which $s_i$ or $\wt s_i$ is equal to 0, in which case we have $\sgn(s_i)
	\cdot \sgn(\wt s_i)=0$). In this
case, by Hoeffding's inequality, the test declares a match with
probability $e^{-\Omega(\ell/\lambda)}$. Let us call this situation a
``spurious match''.

On the other hand, if one of the bits in the $i$-th
block of $\wtt$ was copied over from somewhere in the $i$-th block of
$\w$ for some $i$, we consider the match to be a ``true
match''. Let us describe one
relatively likely way in which this could happen. Suppose that a positive fraction of
the bits in the $i$-th block of $\wtt$ came from the $i$-th block of
$\w$. There is roughly a $e^{-O(\ell/\lambda^2)}$ chance that this
will happen for all $i$ (comparable to the probability that a simple
random walk stays within $[-\lambda, \lambda]$ for $\ell$ steps). In
this case, it is quite likely for a match to occur, because each
$\sgn(s_i)$ and $\sgn(\wt s_i)$ will be positively correlated.

To summarize, the main feature of our test is that it has a true match
rate of at least $e^{-O(\ell/\lambda^2)}$ and a spurious match rate of
at most $e^{-\Omega(\ell/\lambda)}$, which is much lower than the true
match rate as long as $\lambda$ is large enough. In order to make
these statements rigorous, however, the real test we use (as well as
the precise notions of ``true'' or ``spurious'' match) is slightly
more complicated than what is described above. The test is defined
formally in Section \ref{sec:test}.

\subsection{Alignment}

We first remark that our alignment procedure will actually fail for
most (all but about $e^{-O(\log^{1/3} n)}$) of the traces, and we will
disregard these failed traces for purposes of reconstructing the
current bit. Nevertheless, by taking enough total traces (more precisely, by
choosing the constant $M$ in Theorem \ref{thm:main} sufficiently
large), we will still have a sufficient number of successful
alignments to work with.

\begin{figure}
	\centering
	\includegraphics[scale=1]{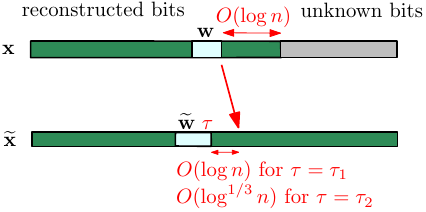}
	\caption{Illustration of the two alignment steps. In each step we fix a substring $\w$ of $\x$ which has distance $O(\log n)$ from the first unknown bit. Assuming the test gives a positive result with some substring $\wt\w$ of the trace, the alignment error (indicated by the lower horizontal arrow) is the difference between the right end-point $\tau$ of $\wt\w$ and the position in $\wt\x$ corresponding to the right end-point of $\w$ (indicated by the vertical red arrow).}
	\label{fig:alignment}
\end{figure}

The alignment $\tau$ is computed in two steps. See Figure \ref{fig:alignment}. We first compute a
preliminary alignment position $\tau_1$ in the output which
corresponds to position $k_1 := k - C \log n$ in the input (where $C$
is some large enough constant). This is done by performing the test
$T_{\ell,\lambda}$ with $(\ell,\lambda) = (C\log^{5/3} n, C^{1/2}
\log^{2/3} n)$. In particular, we declare $\tau_1$ to be the first
index in the output for which
\[ T_{\ell,\lambda}(\x(k_1 - \ell + 1 : k_1), \xt(\tau_1 - \ell + 1 : \tau_1)) = 1. \]
Assuming that such a $\tau_1$ exists, we claim that $\tau_1$ is very
likely to have alignment error of order at most $O(C \log n)$.

The probability of a spurious match is $e^{-\Omega(\ell/\lambda)} =
e^{-\Omega(C^{1/2} \log n)} = n^{-\Omega(C^{1/2})}$, which is
negligible for large enough $C$. This probability is so small that
even taking a union bound over all substrings of length $\ell$, we are
unlikely to find a single length-$\ell$ substring that produces
spurious matches at a high rate.

Meanwhile, the probability of a true match is at least
$e^{-O(\ell/\lambda^2)} \ge e^{-O(\log^{1/3} n)}$. A true match in
this case means that for some $i$ with $1 \le i \le \ell$, the bit in
position $k_1 - \ell + i$ of the input was copied to somewhere very
close to position $\tau_1 - \ell _ i$ in the output (as will be made
precise later).

\begin{figure}
	\centering
	\includegraphics[scale=1]{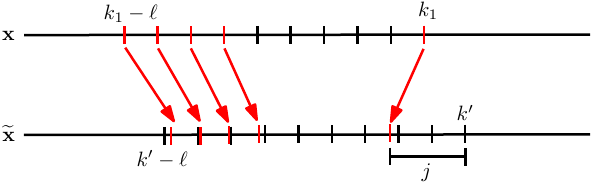}
	\caption{The red arrows indicate to what position a few bits (represented by the  vertical red line segments) in the input string are copied in the trace. The vertical line segments (black or red at the top; black at the bottom) indicate the blocks of length $\lambda$ which are used by the test $T_{\ell,\lambda}$. In the example we may get a positive test with the intervals $\x(k_1-\ell+1:k_1)$ and $\wt\x(k'-\ell+1:k')$ although the alignment error $j$ is rather large, since the number of deletions minus the number of insertions in the latter part of the trace happens to be particularly large, which gives a good overlap between the blocks at the beginning of the string but not at the end.}
	\label{fig:shift}
\end{figure}

However, this does not guarantee that the alignment error of $\tau_1$
is $O(C \log n)$: it could happen that between positions $k_1 - \ell +
i$ and $k_1$, the difference between the number of insertions and
deletions is more than $C \log n$. See Figure \ref{fig:shift} for an illustration of this effect. By Hoeffding's inequality, the
probability this happens is at most $e^{-\Omega\left((C \log n)^2/\ell
  \right)} = e^{-\Omega(C \log^{1/3} n)}$, so that with $C$ large
enough, this misalignment scenario happens rarely compared to the true
match probability described above.

\subsection{Fine alignment and bit statistics}
\label{sec1}

Once we have aligned to within $O(\log n)$, we will perform a second
alignment step to align within roughly $O(\log^{1/3} n)$. This
involves performing the test $T_{\ell,\lambda}$ with $(\ell,\lambda) =
(C^{2/3} \log^{1/3} n, C^{1/12})$. This gives a true match rate of
$e^{-O(\ell/\lambda^2)} = e^{-O(C^{1/2} \log^{1/3} n)}$ and spurious
match rate of $e^{-\Omega(\ell/\lambda)} = e^{-\Omega(C^{7/12}
  \log^{1/3} n)}$, so the true matches dominate the spurious ones.

However, the chance of spurious matches is not low enough to union
bound over all substrings of length $\ell = C^{2/3} \log^{1/3} n$ that
we might align to. In fact, it is actually quite likely that there
will be some ``bad'' substring of length $\ell$ in the input that
produces a high rate of spurious matches (for example, imagine having
two nearby runs of $\ell$ consecutive $0$'s in the input). Getting
around this problem requires some care; the main idea is that although
there may be some bad length-$\ell$ substrings, it is very unlikely
that every single length-$\ell$ substring in an interval of size $\log
n$ is bad, and we use one of the not bad strings to align.

The end result of our two alignment steps is some position $\tau_2$ in
the output which has an alignment error of $O(\log^{1/3} n)$ to some
specified location $k_*$ in the input string $\x$ (i.e., the input position corresponding to
$\tau_2$ is within $O(\log^{1/3} n)$ of $k_*$). Furthermore, $k_*$ is
within the last $O(\log n)$ positions of what we have reconstructed so
far (i.e., $k - k_* = O(\log n)$).

We can then use a variant of the worst case reconstruction algorithm from
\cite{DOS16,NP16} (which was also used in \cite{PZ17}) to reconstruct
$\x(k_*:k_* + C\log n)$ using $e^{O(\log^{1/3} n)}$ (approximate)
traces $\xt(\tau_2:\infty)$, where shifting of $O(\log^{1/3} n)$ can
be tolerated. The algorithm is based on looking at individual bit
statistics, and the key property can be roughly stated as follows:
consider two possibilities for the string $\x$ which match the first
$k$ bits reconstructed so far but disagree on the $(k+1)$-th
bit. Then, there will be a noticeable ($e^{-O(\log^{1/3} n)}$)
difference in the expected value of one of the positions in
$\xt(\tau_2:\infty)$, which allows us to statistically distinguish the
two possibilities. A precise formulation is given in Lemma
\ref{lem:approx-lp}.

\subsection{Implementing the algorithm efficiently}

While we have thus far given a complete description of how to do trace
reconstruction using $e^{O(\log^{1/3} n)}$ traces, there are several
obstacles to making this algorithm run in $n^{1 + o(1)}$ time.

First, during the alignment stage as described so far, we perform our
test $T$ on a sliding window that potentially passes over the whole
output string. This is at least $O(n)$ work needed for reconstructing
even a single bit, which would lead to an overall running time no
better than $O(n^2)$. However, it is quite wasteful to compute our
alignment from scratch each time we reconstruct a new bit. Instead, we
can use previous alignments as a rough guide, allowing us to skip past
all but $n^{o(1)}$ bits of the output string.

Next, to determine the good index $k_2$, we need some way of assessing
whether a given string of length $\log^{1/3} n$ behaves well with our
test $T$ in terms of having a high true match rate and low spurious
match rate. While explicitly calculating these probabilities is not
straightforward, we can estimate them by Monte-Carlo
simulation. Recall that the probabilities in question are of order
$e^{-O(\log^{1/3} n)} = n^{-o(1)}$, so only $n^{o(1)}$ samples are
required to achieve a good enough accuracy.

Finally, in the actual reconstruction step (based on bit statistics as described in the last paragraph of Section \ref{sec1}), a naive implementation
requires us to test every possibility for the first roughly $\log n$
unreconstructed bits. However, as observed in \cite{HMPW08}, the
comparison of bit statistics may be formulated as a linear program,
which can be solved much more efficiently.

\section{Notation for the insertion-deletion channel and Markov properties}
\label{sec:notation}

Let us introduce some general notation and conventions that will be
used for the rest of the paper. To lighten notation, we fix once and
for all two values $q, q' \in [0, 1)$ for the deletion and insertion
  probability, respectively. In order to simplify notation, we will
  further assume throughout that
\begin{equation} \label{eq:q=q'}
  q = q',
\end{equation}
so that the expected length of the output equals the length of the
input. The same arguments carry through in a straightforward way for
general $q, q'$ with appropriate scaling of output lengths. We will
use big-$O$ notation, and all implicit constants in $O(\,\cdot\,)$ and
$\Omega(\,\cdot\,)$ expressions may depend on $q$ and $q'$ but nothing
else.

We will also need to control the relative sizes of various constant
factors. To this end, we introduce a parameter $C$ which will appear
in some of our bounds, which should be thought of as a ``large
constant''. We will ultimately complete our argument by choosing $C$
to be sufficiently large (where the threshold for being large enough
depends only on $q$ and $q'$).

Next, let us introduce some notation relating to strings and their 
traces. Recall that $\cS := \{ 0 , 1 \}^\N$ denotes the space of
infinite sequences of zeroes and ones. Let $\Omega = \cS \times
[0,1]^{\N}$.  We denote the first coordinate function on $\Omega$ by
$\x := (x_0 , x_1 , \ldots)$ and the second by $\omega := (\omega_0 ,
\omega_1 , \ldots)$. Let $U$ be the product uniform measure on
$[0,1]^\N$. If $\rho$ is any measure on $\{ 0 , 1 \}^\N$, let $\P_\rho
:= \rho \times U$. Thus, our previous notation can be expressed as
$\P_\x := \P_{\delta_\x}$ and $\P := \P_\mu$, where $\mu$ is the law
of i.i.d.\ Bernoulli random variables with parameter $1/2$.

We can construct the output $\xt=(\wt x_0,\wt x_1,\dots)$ of the
insertion-deletion channel as a function of $\x$ and $\omega$, where
$\omega$ represents the (random) pattern of insertions and
deletions. The construction proceeds as follows. Temporarily denote $a
:= q(1-q')/(1-qq')$ and $b := q'(1-q)/(1-qq')$. For each $m \in \N$ we
define quantities $s(m) , s'(m) \in \N$, where $s(m)$ (resp.\ $s'(m)$)
represents a position in $\x$ (resp.\ $\xt$) associated with the
randomness of $\omega(m)$. We make the definition by setting $s(0) =
s'(0) = 0$ and proceeding inductively for $m \ge 0$:

\begin{itemize}
	\item If $\omega(m) \in[0,a]$, then define $s(m+1)=s(m)+1$ and
	$s'(m+1)=s'(m)$ (deletion).
	\item If $\omega(m)\in(a,a+b/2]$, then set $s(m+1)=s(m)$, $s'(m+1)=s'(m)+1$,
	and $\wt x_{s'(m)}=0$ (insertion of 0).
	\item If $\omega(m)\in(a+b/2,a+b]$, then set $s(m+1)=s(m)$, $s'(m+1)=s'(m)+1$,
	and $\wt x_{s'(m)}=1$ (insertion of 1).
	\item If $\omega(m)\in(a+b,1]$, then set $s(m+1)=s(m)+1$, $s'(m+1)=s'(m)+1$,
	and $\wt x_{s'(m)}=x_{s(m)}$ (copy).
\end{itemize}

We will now justify briefly why this definition of the deletion-insertion channel is equivalent to the one given in Section \ref{sec:reconstruction}. The channel described in Section \ref{sec:reconstruction} is equivalent to the following: (i) before bit $x_j$ insert $G_j$ uniform and independent bits, where $G_j$ is as before, (ii) delete each of the inserted bits independently with probability $q$, and (iii) delete each of the original bits independently with probability $q$. The combined effect of (i) and (ii) is to insert $\wh G_j$ uniform and independent bits before bit $x_j$, where $\wh G_j$ is a geometric random variable with parameter $b$. From this we conclude equivalence with the channel as described in the bullet points above, because in this channel we insert a geometric number of bits with parameter $b$ between each copy or deletion, and since the fraction of bits in the input string which are deleted is given by $b/(1-a-b)=q$. 

\begin{figure}[ht!]
	\centering
	\framebox{\includegraphics[scale=1]{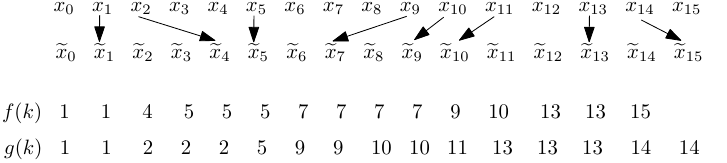}}
	\caption{Illustration of the functions $f$ and $g$.
		The arrows indicate bits which are copied from $\x$ to $\xt$.}
	\label{fig:fg}
\end{figure}

Let us now introduce some notation for corresponding positions in
input strings with positions in their traces. Define
\[ \psi(j) := \inf\{ t\geq 0 \,:\, s(t) = j \}, \qquad \wt{\psi}(j) := \sup\{ t\geq 0 \,:\, s'(t) = j \}. \]
In other words, $\psi(j)$ is the index of the first coordinate in
$\omega$ that decides whether to insert a bit before $x_j$. Similarly,
$\wt{\psi}(j)$ is the index in $\omega$ that determines the value of
$\wt x_j$. Next, define
\[ f(k) = s'(\psi(j)), \qquad g(k') = s(\wt{\psi}(k')). \]
In other words, $f(k)$ is the value of $s'(m)$ at the first time when
$s(m) = k$, and $g(k')$ is the value of $s(m)$ at the first time when
$s'(m) = k'$. See Figure~\ref{fig:fg} for an illustration. Roughly speaking, position $k$ in the
input gets mapped to position $f(k)$ in the output, and position $k'$
in the output was mapped to from position $g(k')$ in the input. In
particular, $g$ is an approximate inverse of $f$. We also define the
function
\begin{equation}
 d(k, k') = \max\left( |f(k) - k'|, |g(k') - k| \right),
 \label{eq:d}
\end{equation}
which measures the failure of position $k$ in the input to correspond
to position $k'$ in the output.

The next proposition records the Markov property that is satisfied
by our insertion-deletion process. It will be convenient to use the
notation $\theta$ to denote the shift operator on bit strings,
i.e., $\theta(\x) = \x(1:\infty)$ and more generally, $\theta^k(\x) =
\x(k:\infty)$.

\begin{proposition} \label{prop:markov-property}
  For any $m \ge 0$, conditioned on the values of $s(m)$ and $s'(m)$,
  the string $\theta^{s'(m)}(\xt)$ has the same law as a trace from
  $\theta^{s(m)}(\x)$ through the insertion-deletion channel. In
  particular, for any $t \ge 0$, the law of $\theta^t(\xt)$ is
  $\P_{\theta^{g(t)}(\x)}$.
\end{proposition}
\begin{proof}
  This is almost immediate from the way we contructed the
  insertion-deletion channel in terms of $\omega$. It is clear that
  the increments of $(s(m), s'(m))$ are i.i.d. Thus, starting $(s(m),
  s'(m))$ from a specified value simply amounts to ignoring the first
  $s(m)$ bits of the input and writing to the output starting from
  position $s'(m)$.
\end{proof}

\section{Clear robust bias test} 
\label{sec:test}
We now give a formal definition of the test $T_{\ell,
  \lambda}$. Recall that the test is designed to answer whether a
substring $\wtt$ of length $\ell$ in a trace is likely to have come
from a substring $\w$ of the same length in the already recovered part
of the input. The test involves subdivision into ``blocks'' of size
approximately $\lambda \leq \sqrt{\ell}$. We remark that the test will
be applied for $(\ell, \lambda)$ on two different scales, namely of
order $(\log^{5/3} n, \log^{2/3} n)$ and $(\log^{1/3} n, 1)$.

Given a string $\w:= \x(k-\ell+1:k)$, let $d := \ceil{\ell / \lambda}$
denote the number of blocks. The right endpoints of the blocks $\{ u_i
\}$ will be given by $u_i := k - \ell + \lceil i \ell / d
\rceil$. Because $\lambda \leq \sqrt{\ell} \leq d$, this definition
makes $\{ (u_{i-1} , u_i] : 1 \leq i \leq d \}$ a partition of $(k -
\ell , k]$ into consecutive intervals of length $\lambda$ or $\lambda
- 1$.

\begin{figure}[h]
	\centering
	\includegraphics[scale=1]{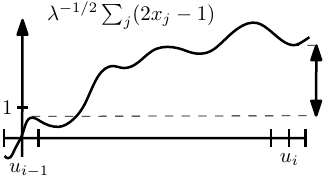}
	\caption{The length of the vertical arrow describes the robust bias
		associated with the block $\x(u_{i-1} + 1 : u_i)$.  The curve
		represents the partial sums $\lambda^{-1/2}\sum_j (2x_j-1)$,
		renormalized to equal 0 at $u_{i-1}$. We say that the robust bias
		is clear if it is at least 1, such as shown in the given example.}
	\label{fig:robust}
\end{figure}

Let us define the {\bf robust bias} of a block $\x(u_{i-1}+1 : u_i)$
to be
\begin{equation} \label{eq:robust-bias}
  \lambda^{-1/2}\inf_{\substack { \text{$t_1,t_2\in\N\,:\,$}
      \\ \text{$|t_1-u_{i-1}|<\lambda/100$}
      \\ \text{$|t_2-u_{i}|<\lambda/100$} \\ }} \Big|\sum_{j=t_1}^{t_2}
  (2x_j-1)\Big|.
\end{equation}
We say that a block has a {\bf clear robust bias} if its robust bias
is at least 1. See Figure~\ref{fig:robust} for an illustration.

For some $\theta\in(0 , 1/10)$ let $\cI \subset \{1 , \dots , d \}$ be
the indices of the $\lceil \theta d \rceil$ blocks for which the
robust bias is largest (with ties resolved in some arbitrary way). By
Donsker's theorem, for $\theta$ sufficiently small and $\lambda$
sufficiently large compared to $\theta$, it holds with high
probability for large $\ell$ that all blocks in $\cI$ have a clear
robust bias. We fix such a choice of $\theta$ as follows: for $B$ a
standard Brownian motion, let
\begin{equation} \label{eq:theta}
  \theta := \frac{1}{10} \P \left[
    \inf_{t_1 \in [0 , 1/50] , t_2 \in[1 , 1 + 1/50]} |B_{t_2} - B_{t_1}| > 1
  \right] > 0 \, .
\end{equation}

For each $i$, define the quantity
\begin{equation} \label{eq37}
  s_i := \sum_{j = u_{i-1}+1}^{u_i} (2x_j - 1)
\end{equation}
which counts the number of 1's minus the number of 0's in the $i$-th
block. Define $\wt s_i$ similarly for a string $\wtt$ of the same length as $\w$. We first define our test using an extra parameter $c \in (0,
1)$. The test is given by
\[ T^c_{\ell, \lambda}(\w, \wtt) =
\begin{cases}
  1 & \text{if } \displaystyle \sum_{i \in \cI} \sgn (s_i) \cdot \sgn(\wt{s}_i) > c |\cI|,\\
  0 & \text{otherwise},
\end{cases} \]
where $|\cI|$ denotes the cardinality of $\cI$.
We will only apply the test for a particular value of $c$, which will
be chosen depending on the insertion/deletion probabilities as
described in the next subsection.

\subsection{Estimates for the test $T$}

\begin{definition}
  We say that a string $\wtt$ of length $\ell$ {\bf has clear robust
    bias at scale $\lambda$} if at least $\lceil \theta d \rceil =
  |\cI|$ of its blocks have clear robust bias.
\end{definition}

First, we formally state our earlier claim that due to our
sufficiently small choice of $\theta$, a random string has clear
robust bias with high probability.

\begin{lemma} \label{lem:robust-bias}
  Let $\w$ be a random string of length $\ell$. Then, $\w$ fails to
  have clear robust bias at scale $\lambda$ with probability at most
  $e^{-\Omega(\ell/\lambda)}$.
\end{lemma}
\begin{proof}
  Within a single block, by Donsker's theorem, the partial sums of the
  bits converge to Brownian motion as $\lambda \rightarrow
  \infty$. Thus, for sufficiently large $\lambda$, the probability
  that this block has clear robust bias is at least
  $2\theta$. Consequently, the probability that the proportion of
  blocks with clear robust bias is less than $\theta$ is exponentially
  small in the number of blocks, i.e., it is of order
  $e^{-\Omega(\ell/\lambda)}$.
\end{proof}

\begin{lemma} \label{lem:T-match-lower}
  Let $\w$ be a string of length $\ell$ which exhibits robust bias at
  scale $\lambda$. Suppose that we take a trace $\wtt$ of $\w$ through
  the insertion-deletion channel. Then, for a small enough constant $c
  > 0$ depending only on the insertion/deletion probabilities, we have
  \[ \P\left( T^c_{\ell,\lambda}(\w, \wtt(0:\ell-1)) = 1 \right) \ge e^{-O(\ell/\lambda^2)}. \]
\end{lemma}

In order to prove Lemma \ref{lem:T-match-lower}, we first define a
condition describing when $\w$ and $\wtt$ are unusually
well aligned. We will then show that $T^c_{\ell,\lambda}$ is
reasonably likely to test positive conditioned on this good alignment.

\begin{definition}
  Let $\w$ be a string of length $\ell$ and $\wtt$ a trace of $\w$
  through the insertion-deletion channel. Following the notation of
  this section, for a block size $\lambda$, we let $d = \lceil
  \ell/\lambda \rceil$, and we let $\{u_i\}_{i = 0}^d$ denote the
  endpoints of the blocks. We say that the trace is {\bf $s$-aligned}
  if for each $0 \le i \le d$, it holds that
  \[ \left| (u_i - u_0) - (f(u_i) - f(u_0)) \right| < s. \]
  If the trace is $s$-aligned, we say that we have {\bf
    $s$-alignment}.
\end{definition}

\begin{proof}[Proof of Lemma \ref{lem:T-match-lower}]
  Let $c_0$ be a small constant to be specified later. We first show
  that the probability of $c_0 \lambda$-alignment is at least
  $e^{-O(\ell/c_0 \lambda^2)}$. For convenience, define the function
  \[ \zeta(x) := (x - u_0) - (f(x) - f(u_0)). \]
  Consider any $\lambda$ consecutive blocks numbered $i+1$ through
  $i+\lambda$. Note that the sequence $\left\{
  \frac{\zeta(k)}{\lambda} \right\}_{k=u_i}^{u_{i+\lambda}}$ is a
  mean-zero simple random walk (since we assume $q=q'$), where the
  distribution of the increments is determined by the insertion and
  deletion probabilities and has finite variance.

  By Donsker's theorem, $\left\{ \frac{\zeta(k)}{\lambda}
  \right\}_{k=u_i}^{u_{i+\lambda}}$ converges to a constant multiple
  of standard Brownian motion as $\lambda \rightarrow \infty$ (where
  time is scaled by $\lambda^2$). For a standard Brownian motion $\{
  W_t \}_{0 \le t \le 1}$, we have that
  \[ \P\Big( W_0, W_1 \in [-c_0/2, c_0/2] \text{ and $|W_t| \le c_0$ for all $0 \le t \le 1$} \Big) = e^{-O(1/c_0)}, \]
  Thus, we have a similar statement for the quantities
  $\frac{\zeta(k)}{\lambda}$, where
  \[ \P\left(
  \frac{|\zeta(u_i)|}{\lambda} \le \frac{c_0}{2}, \;\;
  \frac{|\zeta(u_{i+\lambda})|}{\lambda} \le \frac{c_0}{2}, \;\;
  \text{and} \;\; {\displaystyle \max_{i < j < i + \lambda} \frac{|\zeta(u_{i+\lambda})|}{\lambda} \le c_0} \right) \ge e^{-O(1/c_0)}. \]
  Chaining together these events for each group of $\lambda$
  consecutive blocks (of which there are roughly $\ell/\lambda^2$), we
  see that the overall probability of being $c_0\lambda$-aligned is at
  least $e^{-O(\ell/c_0\lambda^2)}$.

  Our next step is to estimate the probability of a positive test
  conditioned on the trace being $c_0 \lambda$-aligned. We further
  condition on the specific values of the $f(u_i)$. Note that the
  substring $\w(u_i+1:u_{i+1})$ must be transformed into the
  substring $\wtt(f(u_i)+1:f(u_{i+1}) )$. The insertion/deletion
  patterns of these transformations are all independent, and they have
  the same distribution as the insertion-deletion channel applied to
  $\w(u_i:u_{i+1} - 1)$ conditioned on the output trace having length
  $f(u_{i+1}) - f(u_i)$.

  Let $m_i$ denote the number of bits copied to
  $\wtt(f(u_i):f(u_{i+1})-1)$ from $\w(u_i:u_{i+1}-1)$, and note that
  we have $m_i = \Omega(\lambda)$ with probability $1 -
  e^{-\Omega(\lambda)}$. Also, since we assume that our trace is
  $c_0 \lambda$-aligned, only $O(c_0 \lambda)$ bits in
  $\wtt(u_i:u_{i+1}-1)$ were copied from outside of
  $\w(u_i:u_{i+1}-1)$. Finally, note that any bits in
  $\wtt(u_i:u_{i+1}-1)$ which were not copied from $\w$ are
  i.i.d.\ uniformly random.

  Recall that $\wtt$ was assumed to have robust bias at scale
  $\lambda$. Thus, if $i \in \cI$, then with probability $1 -
  e^{-\Omega(\lambda)}$ we will have $\Omega(\sqrt{\lambda})$ bias in
  the sum of bits in $\wtt(u_{i-1}:u_i)$ copied from $\w$. It follows
  that for large enough $\lambda$ and small enough $c_0$, the
  correlation between $\sgn(\wt{s}_i)$ and $\sgn(s_i)$ is
  $\Omega(1)$. This means that (after all our conditioning) as long as
  $c$ is small enough, the test $T^c_{\ell, \lambda}$ will succeed
  except on an event with probability $e^{-\Omega(\ell/\lambda)}$.

  In summary, there is at least a $e^{-O(\ell/\lambda^2)}$ chance to
  have $\lambda$-alignment, and conditioning on this, a positive match
  occurs with high probability. This proves the lemma.
\end{proof}

From now on, we will always run our tests using the value of $c$ as in
Lemma \ref{lem:T-match-lower}. Thus, we will henceforth simply write
$T_{\ell,\lambda}$ instead of $T^c_{\ell,\lambda}$.

\begin{definition} \label{def:mismatch}
  Let $\w = \x(a+1:a + \ell)$ be a substring of the input, and let $\wtt
  = \xt(b+1:b + \ell)$ be a substring of the same length taken from the
  trace. We say that $\w$ and $\wtt$ are \textbf{$s$-mismatched} if
  for any $0 \le i \le \ell$ it holds that $d(a + i, b + i) > s$, where $d$ is defined in \eqref{eq:d}.
\end{definition}

\begin{lemma} \label{lem:T-bad-match-upper}
  Let $\x$ be a random string, and suppose we sample a trace $\xt$
  from the insertion-deletion channel.

  Consider two length-$\ell$ substrings $\w = \x(a+1:a+\ell)$ and $\wtt
  = \xt(b+1:b+\ell)$. Let $\omega^0$ be a realization of the randomness
  of the insertion-deletion channel (i.e., an insertion/deletion
  pattern) for which $\w$ and $\wtt$ are $\lambda$-mismatched. Then,
  \[ \P \Big( T_{\ell,\lambda}(\w, \wtt) = 1 \,\Big|\, \omega = \omega^0 \Big) \le e^{-\Omega(\ell/\lambda)}. \]
\end{lemma}
\begin{proof}
  Note that after conditioning on $\omega$, the remaining randomness
  is to sample the values of the bits. We do this incrementally block
  by block (where the blocks are of size $\lambda$). Suppose that we
  have already sampled the first $i - 1$ blocks of $\w$ and $\wtt$;
  let us consider the conditional distribution on bits in the $i$-th
  blocks of $\w$ and $\wtt$.

  Note that there are only two ways dependency can occur between bits
  in the $i$-th blocks of $\w$ and $\wtt$ and the already sampled
  bits: either a bit in the $i$-th block of $\wtt$ came from one of
  the first $i - 1$ blocks of $\w$ or a bit in one of the first $i -
  1$ blocks of $\wtt$ came from a bit in the $i$-th block of
  $\w$. Moreover, these two scenarios are mutually exclusive. It
  follows that the $i$-th block of at least one of $\w$ and $\wtt$ is
  completely independent of the already sampled bits.

  Suppose for example that the $i$-th block of $\w$ is independent of
  the already sampled bits. Since $\w$ and $\wtt$ are
  $\lambda$-mismatched, the $i$-th blocks of $\w$ and $\wtt$ are also
  independent of each other. It follows that we may first sample the
  $i$-th block of $\wtt$, and conditioned on that, $\sgn(s_i)$ is
  equally likely to be $1$ or $-1$. A similar argument holds in the
  other case, where the $i$-th block of $\wtt$ is independent of the
  first $i - 1$ blocks of $\w$ and $\wtt$.

  Either way, the end result is that $\sgn(s_i) \cdot \sgn(\wt{s}_i)$
  is a uniform random sign independent of the previous blocks. It
  follows that the quantity
  \[ \sum_{i \in \cI} \sgn(s_i) \cdot \sgn(\wt{s}_i) \]
  has the law of a sum of $|\cI|$ independent random signs, so by
  Hoeffding's inequality, the chance that it exceeds the threshold
  $c|\cI|$ required for our test $T_{\ell,\lambda}$ is
  $e^{-\Omega(|\cI|)} = e^{-\Omega(\ell/\lambda)}$, as desired.
\end{proof}

\section{The ``good'' set of strings}

\begin{definition} \label{def:spurious}
  Let $\ell$ and $\lambda$ be given positive integers. Let $\x$ be an
  input string, let $I$ be an interval of length $\ell$, and write $\w
  = \x(I)$. Let $J$ be another interval (often we will have $I
  \subseteq J$).

  Suppose we take a trace of $\x$ through the insertion-deletion
  channel. We say that an \textbf{$(I, J)$-spurious match} occurs if
  for some substring $\wtt=\xt(i_1:i_2)$ of the output such that $g([i_1,i_2]) \subseteq
  J$, we have $T_{\ell,\lambda}(\w, \wtt) = 1$, but $\w$ and $\wtt$
  are $\lambda$-mismatched. We use
  \[ \cQ_{T_{\ell, \lambda}}(I, J) \]
  to denote the event that an $(I, J)$-spurious match occurs.
\end{definition}

\begin{lemma} \label{lem:spurious-bound}
  Let $\ell$ and $\lambda$ be given positive integers. Let $I$ be an
  interval of length $\ell$, and let $J \supseteq I$ be an interval
  containing $I$. Suppose we have an input string $\x$ all of whose
  bits are determined except those in $J$, which are drawn
  i.i.d.\ uniformly. Then, letting $|J|$ denote the length of $J$,
  \[ \P(\cQ_{T_{\ell, \lambda}}(I, J)) \le |J| e^{-\Omega(\ell/\lambda)} + e^{-\Omega(|J|)}. \]
\end{lemma}
\begin{proof}
  We will first condition on an insertion/deletion pattern $\omega$
  from the insertion-deletion channel. For each interval $I'$ of
  length $\ell$, if $\omega^0$ is an insertion/deletion pattern which
  causes $I'$ and $I$ to be $\lambda$-mismatched, then we know by
  Lemma \ref{lem:T-bad-match-upper} that
  \begin{equation} \label{eq:spurious-bound-single-pair}
    \P(T_{\ell,\lambda}(\x(I), \xt(I')) = 1 \mid \omega = \omega^0) \le e^{-\Omega(\ell/\lambda)}.
  \end{equation}
  Let $J'$ denote the minimal interval containing $f(J)$, and note
  that $J'$ is a function of the insertion/deletion pattern
  $\omega$. Then, conditioning on $J'$ and taking a union bound over
  all possible $I' \subseteq J'$ in
  \eqref{eq:spurious-bound-single-pair} gives
  \[ \P(\cQ_{\ell,\lambda}(I, J) \mid J') \le |J'| e^{-\Omega(\ell/\lambda)}. \]
  Since we also have $\P(|J'| \ge 2|J|) \le e^{-\Omega(|J|)}$, this
  yields our final bound
  \[ \P(\cQ_{T_{\ell, \lambda}}(I, J)) \le |J| e^{-\Omega(\ell/\lambda)} + e^{-\Omega(|J|)}. \]
\end{proof}

\begin{lemma} \label{lem:disjoint-bound}
  Let $\ell$ and $\lambda$ be given positive integers. Let $I$ be an
  interval of length $\ell$, and let $J$ be another disjoint interval
  whose distance from $I$ is at least $|J|$. Suppose we have an input
  string $\x$ all of whose bits are determined except those in $I$,
  which are drawn i.i.d.\ uniformly. Then,
  \[ \P(\cQ_{T_{\ell, \lambda}}(I, J)) \le |J| e^{-\Omega(\ell/\lambda)} + e^{-\Omega(|J|)}. \]
\end{lemma}
\begin{proof}
  Let $J'$ be the minimal interval containing $f(J)$, and define the
  event
  \[ E = \left\{ \text{$|J'| \le 2|J|$ and not all bits between $I$ and $J$ were deleted} \right\}. \]
  Note that $E$ is measurable with respect to the $\sigma$-field
  generated by $\omega$ (i.e., the randomness of the insertion-deletion
  channel), and $\P(E) \ge 1 - e^{-\Omega(|J|)}$.

  Meanwhile, conditioned on $E$, consider any subinterval $I'
  \subseteq J'$ of length $\ell$. None of the bits in $I'$ come from
  $I$, so the random bits in $I$ are independent of the bits in
  $I'$. Consequently, we have
  \[ \P(T_{\ell, \lambda}(\x(I), \xt(I')) = 1 \mid E) \le e^{-\Omega(\ell/\lambda)}. \]
  Taking a union bound over possible choices of $I'$ given that $E$
  occurs, we conclude that
  \[ \P(\cQ_{T_{\ell, \lambda}}(I, J)) \le |J| e^{-\Omega(\ell/\lambda)} + \P(E^c) \le |J| e^{-\Omega(\ell/\lambda)} + e^{-\Omega(|J|)}. \]
\end{proof}

\subsection{Coarsely well-behaved strings}

\begin{definition} \label{def:coarsely-behaved}
  Let $\x$ be a string of length at least $n$. Let $\ell = C \log^{5/3} n$ and
  $\lambda = C^{1/2} \log^{2/3} n$. We say that $\x$ is
  \textbf{coarsely well-behaved} if for each interval $I \subseteq [0,
    n]$ of length $\ell$, it holds that $\x(I)$ has robust bias at scale
  $\lambda$ and
  \[ \P_\x \left( \cQ_{\ell,\lambda}(I, [0, n]) \right) \le n^{-2}. \]
\end{definition}

\begin{lemma} \label{lem:coarsely-behaved}
  Let $\x$ be a random string. Then, $\x$ is coarsely well-behaved
  with probability at least $1 - n^{-2}$.
\end{lemma}
\begin{proof}
  Consider first a particular interval $I = (a, a + \ell]$ of 
  length
  $\ell$. By Lemma \ref{lem:robust-bias}, we know that $\x(I)$ has
  robust bias at scale $\lambda$ with probability at least
  \[ 1 - e^{-\Omega\left( \ell/\lambda \right)} = 1 - e^{-\Omega(C^{1/2} \log n)} \ge 1 - n^{-4} \]
  for large enough $C$.

  We also have from Lemma \ref{lem:spurious-bound} that
  \[ \P(\cQ_{\ell,\lambda}(I, [0, n])) \le e^{-\Omega(C^{1/2} \log n)}. \]
  Thus, for large enough $C$, we can ensure by Markov's inequality that
  \[ \P\left( \P_{\x}(\cQ_{\ell,\lambda}(I, [0, n])) \ge n^{-2} \right) \le n^{-4}. \]
  Taking a union bound over at most $n$ possible values of $I$
  completes the proof.
\end{proof}

\subsection{Finely well-behaved strings}
Recall Lemmas \ref{lem:spurious-bound} and \ref{lem:disjoint-bound}. Let $c_0>0$ be such that the terms $\Omega(\ell/\lambda)$ and $\Omega(|J|)$ in these lemmas could have been replaced by $10c_0\ell/\lambda$ and $10c_0|J|$, respectively.
\begin{definition} \label{def:finely-behaved}
  Let $\x$ be a string of length $n$, and let $\ell = C^{2/3}
  \log^{1/3} n$ and $\lambda = C^{1/12}$. We say that $\x$ is
  \textbf{finely well-behaved} if for each interval $J := [a, a + C
    \log n] \subseteq [0, n]$ of length $C \log n$, there exists a
  subinterval
  \[ I \subset [a + \tfrac{1}{3}C \log n, a + \tfrac{2}{3} C \log n] \]
  of size $\ell$ such that $\x(I)$ has robust bias at scale $\lambda$
  and
  \[ 
  \P_{\x}\left( \cQ_{\ell,\lambda}(I, J) \right) \le e^{- c_0C^{7/12} \log^{1/3} n}. 
  \]
\end{definition}

\begin{lemma}
  Let $\x$ be a random string. Then, $\x$ is finely well-behaved with
  probability at least $1 - n^{-2}$.
\end{lemma}

\begin{proof}
Throughout the proof the implicit constants in $\Omega(\cdot)$ and $O(\cdot)$ may depend on $c_0$ but not on $C$.
  Fix a particular interval $J = [a, a + C \log n] \subseteq [0, n]$
  of length $C \log n$, and let $\ell = C^{2/3} \log^{1/3} n$ and
  $\lambda = C^{1/12}$ be as in Definition
  \ref{def:finely-behaved}. Consider $m := \tfrac{1}{3} C^{1/3}
  \log^{2/3} n$ disjoint length-$\ell$ subintervals
  \[ I_1, I_2, \ldots, I_m \subset [a + \tfrac{1}{3}C \log n, a + \tfrac{2}{3} C \log n]. \]

  For a given realization of $\x$, we say that $I_i$ is \emph{bad} if
  either it does not have clear robust bias at scale $\lambda$ or it
  holds that
  \begin{equation} \label{eq:bad-pattern}
    \P_{\x}\left( \cQ_{\ell,\lambda}(I_i, J) \right) \ge e^{-c_0C^{7/12} \log^{1/3} n}.
  \end{equation}
  Let $\ell' := C^{1/6} \ell = C^{5/6} \log^{1/3} n$ and define the
  event
  \[ H = \left\{ \parbox{19em}{
    \centering there exists some $t$ such that \\
    $g(t), g(t+\ell) \in I$ and $|g(t) - g(t + \ell)| \ge \ell'$
  } \right\}, \]
  which roughly says that a substring of length $\ell'$ had so many
  deletions that only $\ell$ or fewer bits were left in the output. We
  have that $\P(H) \le e^{-\Omega(\ell')}$.

  As long as $H$ does not occur, then any spurious match counted in
  \eqref{eq:bad-pattern} must have come from an interval $J'$ of
  length $\ell'$, i.e., we have
  \[ \cQ_{\ell,\lambda}(I_i, J) \subseteq H \cup \left( \bigcup_{\substack{\text{$J' \subset J$ an interval} \\ \text{of length $\ell'$}}} \cQ_{\ell,\lambda}(I_i, J') \right). \]
  Thus, by the pigeonhole principle, if $I_i$ is bad due to
  \eqref{eq:bad-pattern}, then there must be some interval $J'_i$ of
  length $\ell'$ for which
  \begin{equation} \label{eq:bad-pair}
    \P_{\x}\left( \cQ_{\ell,\lambda}(I_i, J'_i) \right) \ge e^{-c_0C^{7/12} \log^{1/3} n}.
  \end{equation}
  We say that $(I_i, J'_i)$ is a \emph{bad pair} if either $I_i$ does
  not have clear robust bias at scale $\lambda$, or
  \eqref{eq:bad-pair} holds. In particular, the above discussion shows
  that if $I_i$ is bad, then there is some interval $J'_i$ of length
  $\ell'$ for which $(I_i, J'_i)$ is a bad pair.

  Suppose for the sake of contradiction that with probability at least
  $n^{-2}$ in the randomness of $\x$, all the $I_i$ are bad (i.e., $\x$
  is not finely well-behaved). Each $I_i$ is part of a bad pair $(I_i,
  J'_i)$, and note that there are at most $(C \log n)^m = n^{o(1)}$
  possible values for $(J'_1, \ldots , J'_m)$. Thus, by the pigeonhole
  principle, it must hold for some specific choice of $(J'_1, \ldots ,
  J'_m)$ that
  \begin{equation} \label{eq:bad-pair-contradiction}
    \P\left( \text{$(I_i, J'_i)$ is a bad pair for $i = 1, 2, \ldots ,
      m$} \right) \ge n^{-3}.
  \end{equation}
  Fixing this choice of $(J'_1, \ldots , J'_m)$, we will derive a
  contradiction.

  Let $r := 0.01 C^{1/6} \log^{2/3} n$. To carry out the analysis, we
  inductively define a sequence $i_1, i_2, \ldots , i_{r}$ as
  follows. We take $i_1 = 1$, and for $k \ge 1$, let
  \[ N_k = \bigcup_{j = 1}^k (I_{i_j} \cup J'_{i_j}). \]
  Then, choose $i_{k+1}$ so that $I_{i_{k+1}}$ is distance at least
  $2\ell'$ from $N_k$. Note that the $2\ell'$-neighborhood of $N_k$
  	intersects at most $2k \cdot \lceil 5\ell'/\ell\rceil \leq 12 C^{1/6} k$ of the $I_i$, so such a choice is always possible as long as $k \le r$.

  Let $\cG_k$ be the $\sigma$-field generated by the bits of $\x$
  whose positions are in $N_k$, and let $E_k$ denote the event that
  $(I_{i_k}, J_{i_k})$ is a bad pair. Note that $E_k$ is measurable
  with respect to $\cG_k$.

  First of all, note that whether $I_{i_k}$ has clear robust bias at
  scale $\lambda$ is independent of $\cG_{k-1}$, so by Lemma
  \ref{lem:robust-bias}, we have
  \[ \P\left(\parbox{10em}{\centering
    $I_{i_k}$ does not have clear \\ robust bias at scale $\lambda$
  } \,\middle|\, \cG_{k-1} \right) \le e^{-\Omega(\ell/\lambda)} = e^{-\Omega(C^{7/12} \log^{1/3} n)}. \]

  Next, we will estimate $\P\left( \cQ_{\ell,\lambda}(I_{i_k},
  J'_{i_k}) \,\middle|\, \cG_{k-1} \right)$. Suppose first that
  $I_{i_k}$ and $J'_{i_k}$ are disjoint and distance at least $\ell'$
  apart. Then, by Lemma \ref{lem:disjoint-bound}, we have
  \[ \P \left( \cQ_{\ell,\lambda}(I_{i_k}, J'_{i_k}) \,\middle|\, \cG_{k-1} \right) \le \ell' e^{-10c_0\ell/\lambda} + e^{-10c_0\ell'} \leq e^{-9c_0C^{7/12} \log^{1/3} n}. \]
  If instead $I_{i_k}$ and $J'_{i_k}$ are within distance $\ell'$ of
  each other, then let $J$ be the interval formed by extending
  $I_{i_k}$ on both sides by $2\ell'$, so that $J'_{i_k} \subset
  J$. By our construction, it is also guaranteed that $J$ is disjoint
  from $N_{k-1}$, and so when conditioning on $\cG_{k-1}$, none of its
  bits have been determined yet. Then, we may apply Lemma
  \ref{lem:spurious-bound} to obtain
  \begin{align*}
    \P \left( \cQ_{\ell,\lambda}(I_{i_k}, J'_{i_k}) \,\middle|\, \cG_{k-1} \right) &\le \P \left( \cQ_{\ell,\lambda}(I, J) \,\middle|\, \cG_{k-1} \right) \\
    &\le 3\ell' e^{-10c_0\ell/\lambda } + e^{-10c_0\ell'} = e^{-9c_0C^{7/12} \log^{1/3} n}.
  \end{align*}
  Thus, the above bound holds in either case, and by Markov's
  inequality, this implies
  \[ \P\left( \P_{\x} \left( \cQ_{\ell,\lambda}(I_{i_k}, J'_{i_k}) \right) \ge e^{-c_0C^{7/12} \log^{1/3} n} \,\middle|\, \cG_{k-1} \right) \le e^{-\Omega(C^{7/12} \log^{1/3} n)}. \]
  It follows that
  \begin{align*}
    \P(E_k \mid \cG_{k-1}) &\le \P\left(\parbox{10em}{\centering
    $I_{i_k}$ does not have clear \\ robust bias at scale $\lambda$
    } \,\middle|\, \cG_{k-1} \right) \\
    &\hphantom{==} + \P\left( \P_{\x} \left( \cQ_{\ell,\lambda}(I_{i_k}, J'_{i_k}) \right) \ge e^{-c_0 C^{7/12} \log^{1/3} n} \,\middle|\, \cG_{k-1} \right) \\
    &\le e^{-\Omega(C^{7/12} \log^{1/3} n)}.
  \end{align*}
  Iterating this over $k = 1, \ldots , r$, we finally obtain
  \[ \P\left( E_1 \cap E_2 \cap \cdots \cap E_r \right) \le e^{-\Omega(C^{7/12} \log^{1/3} n) \cdot r} = e^{-\Omega(C^{3/4} \log n)}, \]
  which is smaller than $n^{-3}$ for large enough $C$. This gives our
  desired contradiction of \eqref{eq:bad-pair-contradiction},
  completing the proof.
\end{proof}

\section{Alignment rules}

For the rest of the paper, let $\bad$ denote the set of strings that
either fail to be coarsely well-behaved or finely
well-behaved. 
\begin{lemma} \label{lem:tau_1}
  Let $k\geq\ell$ be a given positive integer, and let $\ell = C \log^{5/3} n$
  and $\lambda = C^{1/2} \log^{2/3} n$ as in Definition
  \ref{def:coarsely-behaved}. Define
\[ \tau^k_1 = \tau^k_1(\xt) := \inf \left\{ k'\in[\ell, 2n] : T_{\ell,\lambda}(\x(k-\ell+1:k), \xt(k'-\ell+1:k')) = 1 \right\}, \]
  where we set $\tau^k_1 = \infty$ if no such $k'$ exists. If $\x
  \not\in \bad$, then
  \[ \P(\tau^k_1 < \infty,\; d(k, \tau^k_1) > \tfrac{1}{10} C \log n) \le e^{-\Omega(C \log^{1/3} n)} \]
  and
  \[ \P(\tau^k_1 < \infty,\; d(k, \tau^k_1) \le \tfrac{1}{10} C \log n) \ge e^{-O(\log^{1/3} n)}. \]
\end{lemma}
\begin{proof}
  Let $I = (k - \ell, k]$, and define the events
  \[ E = \{ \tau^k_1 < \infty \}, \qquad E' = \{ \tau^k_1 < \infty,\; d(k, \tau^k_1) > \tfrac{1}{10} C \log n \}, \]
  \[ F = \left\{ \parbox{22em}{\centering
    the difference between the number of inserted and \\
    deleted bits among those in $I$ is at least $\frac{1}{20} C \log n$
  } \right\}. \]
  Note that by Hoeffding's inequality, we have
  \[ \P_{\x}(F) = e^{-\Omega\left( \frac{(C \log n)^2}{C \log^{5/3} n} \right)} = e^{-\Omega(C \log^{1/3} n)}. \]

  To show the first inequality, which amounts to bounding
  $\P_{\x}(E')$, suppose that $E$ holds but $\cQ_{\ell,\lambda}(I,
  [0, n])$ does not. Then, it must be that the matched string
  $\xt(\tau^k_1 - \ell+1:\tau^k_1)$ is not
  $\lambda$-mismatched. However, if additionally $d(k, \tau^k_1) > C
  \log n$, then $F$ must hold. Thus,
  \begin{align*}
    \P_{\x}(E') &\le \P_{\x}(\cQ_{\ell,\lambda}(I, [0, n])) + \P_{\x}(E' \setminus \cQ_{\ell,\lambda}(I, [0, n])) \\
    &\le n^{-2} + e^{-\Omega(C \log^{1/3} n)} = e^{-\Omega(C \log^{1/3} n)},
  \end{align*}
  establishing the first inequality.

  Next, note that by Lemma \ref{lem:T-match-lower}, a positive match
  will be found (i.e., $E$ will hold) with probability at least
  $e^{-O(\ell/\lambda^2)} = e^{-O( \log^{1/3} n)}$. Thus,
  \[ \P_{\x}(E \setminus E_1) \ge e^{-O( \log^{1/3} n)} - e^{-\Omega(C \log^{1/3} n)} = e^{-O( \log^{1/3} n)}, \]
  establishing the second inequality.
\end{proof}

\begin{lemma} \label{lem:tau_2}
  Let $\x \not\in \bad$ be a string, and let $k$ be a positive
  integer. Suppose that we know $\x(0:k)$, and let $\varphi(t) = t
  e^{\frac{t}{C^{2/3} \log^{2/3} n}}$. Let $\cF_j$ denote the $\sigma$-algebra generated by $\x(0:k)$ and $\xt(0:j)$. Then, there is a position $k_*
  \in [k - C \log n, k]$ and a stopping time $\tau^k_2(\xt)$ for $\cF_j$ for which
  the following properties hold: defining the event
  \[ F^k = F^k_{\x} := \{ \tau^k_2(\xt) < \infty,\; g(\tau^k_2(\xt)) \in [k - C\log n, k] \}, \]
  we have
  \begin{enumerate}
  \item $\P\Big[ \{ \tau^k_2(\xt) < \infty \} \setminus F^k \Big] \le n^{-2}$,
  \item $\P\Big[ F^k \Big] \ge e^{-O(C^{1/2} \log^{1/3} n)}$,
  \item $\E\Big[ \varphi \big(|g(\tau^k_2(\xt)) - k_*| \big) \,\Big|\, F^k \Big] \le O(C^{2/3} \log^{1/3} n)$.
  \end{enumerate}
\end{lemma}

\begin{remark}
  In fact, the proof below implies that the same result holds if
  instead of requiring $\x \not\in \bad$, we only require that its
  first $k$ bits match some string not in $\bad$.
\end{remark}

\begin{proof}
  Let $\ell = C^{2/3} \log^{1/3} n$ and $\lambda = C^{1/12}$, and let
  $I \subset [k - \frac{2}{3} C \log n, k - \frac{1}{3} C \log n]$ be
  the interval guaranteed by Definition \ref{def:finely-behaved}
  (since $\x$ is finely well-behaved). Let $\tau^k_1$ be as in Lemma
  \ref{lem:tau_1}, and consider the interval
  \[ I' = [\tau^k_1(\xt) + \tfrac{1}{6} C \log n, \tau^k_1(\xt) + \tfrac{5}{6} C \log n]. \]
  We define
  \[ \tau^k_2(\xt) := \inf \left\{ k' : [k' - \ell, k'] \subset I' \text{ and } T(\x(I), \xt(k' - \ell : k')) = 1\right\}, \]
  where as usual we set $\tau^k_2(\xt) = \infty$ if no such $k'$
  exists (or if $\tau^k_1(\xt) = \infty$). Our choice of $k_*$ is then
  the right endpoint of $I$.

  To lighten notation, in the rest of the proof we write $\tau^k_1 =
  \tau^k_1(\xt)$ and $\tau^k_2 = \tau^k_2(\xt)$. In several of our
  calculations, it will be convenient to exclude the event
  \[ E := \Big\{ \text{$|g(\tau^k_1 + t) - (k - C \log n + t)| > \tfrac{1}{8}C \log n$ for some $0 \le t \le C \log n$} \Big\}. \]
  We first show that $E$ is a rare event. Recall that by the
  properties of $\tau^k_1$ established in Lemma \ref{lem:tau_1}, we
  have that
  \[ \P\Big( |g(\tau^k_1) - k + C \log n| \ge \tfrac{1}{10} C \log n \Big) \le n^{-2}. \]
  On the other hand, if $|g(\tau^k_1) - k + C \log n| \le
  \tfrac{1}{10} C \log n$, then in order for $E$ to occur, among
  the bits with input positions between $k - \tfrac{11}{10} C \log n$
  and $k$, the difference between the number of deletions and
  insertions must have been at least $\Omega(C \log n)$. This occurs
  with probability at most $n^{-\Omega(C)} \le n^{-2}$ for large
  enough $C$. Thus, we see that $\P(E) \le O(n^{-2})$.

  Let us now establish the properties stated in the lemma. The first
  property immediately follows from our bound on $\P(E)$, since
  whenever $\tau^k_2 < \infty$, we have
  \[ g(\tau^k_1) + \tfrac{1}{6}C \log n \le g(\tau^k_2) \le g(\tau^k_1) + \tfrac{5}{6}C \log n. \]
  Thus, we have $\{ \tau^k_2 < \infty \} \setminus F^k \subseteq E$, and so
  \[ \P(\{ \tau^k_2 < \infty \} \setminus F^k) \le \P(E) \le O(n^{-2}). \]

  Recall from Lemma \ref{lem:T-match-lower} that with probability at
  least $e^{-O(\ell/\lambda^2)} = e^{-O(C^{1/2} \log^{1/3} n)}$, the
  bits coming from $I$ will form a positive match for the test
  $T_{\ell,\lambda}$. Outside of the event $E$, this match will be
  detected by our procedure, and so
  \[ \P(\tau^k_2 < \infty) \ge e^{-O(C^{1/2} \log^{1/3} n)} - \P(E) = e^{-O(C^{1/2} \log^{1/3} n)}. \]
  Subtracting the bound from the first property yields the second
  property.

  For the last property, let us estimate the probability
  \[ \P\Big[ F^k \cap \{ |g(\tau^k_2) - k_*| \ge 2\ell \} \Big]. \]
  There are two possible cases to consider:
  \begin{enumerate}
  \item A spurious match event $\cQ_{\ell,\lambda}(I, [k - C\log n,
    k])$ may occur.  
  \item If there is no spurious match event but $F^k$ holds, then it
    means there is some $t \in [\tau^k_2 - \ell, \tau^k_2]$ for which
    $g(t) \in [k_* - \ell, k_*]$. Then, the only way to have
    $|g(\tau^k_2) - k_*| \ge 2\ell$ is if there were at least $\ell$
    more deletions than insertions in the length-$2\ell$ input
    interval $[k_*, k_* + 2\ell]$.
  \end{enumerate}
  By our choice of the interval $I$, we can estimate the probability
  of the first scenario by
  \[ \P(\cQ_{\ell,\lambda}(I, [k - C\log n, k])) \le e^{-\Omega(C^{7/12} \log^{1/3} n)}, \]
  and the last scenario has probability $e^{-\Omega(\ell)} =
  e^{-\Omega(C^{2/3} \log^{1/3} n)}$. Thus, the overall probability is
  \[ \P\Big[ F^k \cap \{ |g(\tau^k_2) - k_*| \ge 2\ell \} \Big] \le e^{-\Omega(C^{7/12} \log^{1/3} n)}. \]
  and so
  \[ \frac{
    \P\Big[ F^k \cap \{ |g(\tau^k_2) - k_*| \ge 2\ell \} \Big] }{
    \P\Big[ F^k \cap \{ |g(\tau^k_2) - k_*| \le 2\ell \} \Big] } \le
  e^{-\Omega(C^{1/2}\log^{1/3} n)}. \]
  To calculate the relevant expectation, we can divide into cases
  depending on whether $|g(\tau^k_2) - k_*| \le 2\ell$ or not (note that on
  the event $F^k$, we always have $|g(\tau^k_2) - k_*| \le C \log
  n$). This yields
  \begin{align*}
    \E\Big[ \varphi(|g(\tau^k_2) - k_*|) \,\Big|\, F^k \Big] &\le \varphi(2\ell) \cdot \P\Big[ |g(\tau^k_2) - k_*| \le 2 \ell \,\Big|\, F^k \Big] + \\
    &\hphantom{==} \varphi(C \log n) \cdot \E\Big[ |g(\tau^k_2) - k_*| \ge 2 \ell \,\Big|\, F^k \Big] \\
    &= \varphi(2\ell) + e^{-\Omega(C^{1/2}\log^{1/3} n)} \cdot \varphi(C \log n) \\
    &= O(\ell) = O(C^{2/3} \log^{1/3} n),
  \end{align*}
  as desired.
\end{proof}

\section{Reconstruction from approximately aligned strings}\label{sec:bit-test}
Recall from Proposition \ref{prop:markov-property} that we can use
$\xt \left(\tau^{k}_2(\xt) : \infty \right)$ as a trace of
$\x(g(\tau^{k}_2(\xt)):\infty)$. If we had exactly $g(\tau^{k}_2(\xt))
= k$, then the problem would be reduced to worst case reconstruction
of $\x(k:\infty)$. This section adapts the methods of
\cite{DOS16,NP16} to handle imperfect alignment using a similar
approach as \cite{PZ17}. We start with the following definition.

\begin{definition} \label{def:v}
  Let $\tau^k_2$ be as in Lemma \ref{lem:tau_2}. For a bit string $\x$
  and positive integer $k$, let
  \[ V(\xt) := \xt \left(\tau^{k}_2(\xt) : \infty \right). \]
  (We will only be concerned with $V(\xt)$ in cases where
  $\tau^k_2(\xt) < \infty$.) For any string $\w$, let $\x_\w$ denote the concatenation $\x(0:k) :
  \w$, and consider a trace $\xt_{\w}$ drawn from the
  insertion-deletion channel applied to $\x_\w$. Then, define
  \[ 
  v(\w) := \E_{\x_\w} \left[ V(\xt_{\w}) \,\middle|\, F^k_{\x_{\w}} \right]. 
  \]
  Note that $v$ is a linear function of $\w$, so we may extend this
  definition to any $\w \in [0,1]^\N$. Finally, for $m\in\N$ and $\a \in [0,1]^\N$ set
  \begin{equation}\label{eq:maxm}
   \Linfm{\a} := \max_{j \le m} |a_j|. 
  \end{equation}
\end{definition}

Our reconstruction is based on the following lemma.

\begin{lemma} \label{lem:approx-lp}
  Fix a string $\x \not\in \bad$ and consider a constant $C >
  0$. Suppose that $\w \in [0,1]^\N$ satisfies
  \begin{equation} \label{eq:v-bound-hypothesis}
    |v_j(\w) - v_j(\x(k+1:\infty))| \le e^{-\Omega(C \log^{1/3} n)}
  \end{equation}
  for all $j \le C^2 \log n$. Then, for large enough $C$, we must have
  \[ |w_0 - x_{k+1}| < \frac{1}{2}. \]
  Moreover, if $\w = \x(k+1:k+2C^2 \log n)$, then
  \begin{equation} \label{eq:v-nearly-optimal}
    \|v(\w) - v(\x(k+1:\infty))\|_{C^2\log n,\infty} \le O(n^{-2}).
  \end{equation}
\end{lemma}

The core of the proof of Lemma \ref{lem:approx-lp} is contained in the
following lemma about bit statistics for randomly shifted strings.

\begin{lemma} \label{lem:distinguishing-with-shifts}
  Let $n$ be a positive integer, and let $\a \in [-1,1]^{\N}$ be a
  sequence of real numbers for which $a_i = 0$ for $i \le n$ but for
  which $|a_{n+1}| \ge \frac{1}{2}$. Let $S$ be a random variable
  taking integer values between $0$ and $n - 1$.

  Let $\varphi(t) = t e^{t/n^{2/3}}$. Then, there exists an index $j =
  O(n)$ such that for a trace $\at$ from the shifted sequence
  $\theta^S(\a)$, we have
  \[ \Big| \E\left[ \E_{\theta^S(\a)}(\wt{a}_j) \right] \Big| \ge \exp\left( -O\left(n^{1/3} + \min_{0 \le t \le n} \E\Big[ \varphi(|S - t|)\Big] \right) \right). \]
\end{lemma}

Let us first deduce Lemma \ref{lem:approx-lp} from Lemma
\ref{lem:distinguishing-with-shifts}.

\begin{proof}[Proof of Lemma \ref{lem:approx-lp}]
  For the first claim, note that $v(\w)$ and $v(\x(k+1):\infty)$ are
  calculated by taking expectations involving traces from $\x_\w$ and
  $\x$. Let $\xt_\w$ and $\xt$ denote these traces, and for purposes
  of our analysis, we may suppose that these traces were sampled using
  the same insertion/deletion pattern $\omega$. In this case, the two
  events $F^k_{\x}$ and $F^k_{\x_\w}$ coincide, because they involve
  only the first $k$ bits of $\x$ and $\x_\w$, which are
  identical. Thus, we will subsequently use $F^k$ to refer to either
  $F^k_{\x}$ or $F^k_{\x_\w}$.

  Note also that when $F^k$ holds, we have $\tau^k_2(\xt) =
  \tau^k_2(\xt_\w)$ and $g(\tau^k_2(\xt)) =
  g(\tau^k_2(\xt_\w))$. Accordingly, we will write $\tau =
  \tau^k_2(\xt) = \tau^k_2(\xt_\w)$ when conditioning on $F^k$. Let
  $S$ be a random variable with the same distribution as $g(\tau) - k
  + C\log n$ conditioned on $F^k$. Then, by Proposition
  \ref{prop:markov-property}, we have
  \begin{align}
    v_j(\w) - v_j(\x(k+1:\infty)) &= \E \left[ \wt{x}_{\w, \tau + j + 1} - \wt{x}_{\tau + j + 1} \mid F^k \right] \nonumber \\
    &= \E_{\theta^{k - C \log n + S}(\x_\w)}[\wt{x}_{j + 1}] - \E_{\theta^{k - C \log n + S}(\x)}[\wt{x}_{j + 1}]. \label{eq:lp-approx1}
  \end{align}

  We are now in a position to apply Lemma
  \ref{lem:distinguishing-with-shifts}. Take
  \[ \a = \x_\w(k - C\log n:\infty) - \x(k - C\log n:\infty). \]
  Note that $a_i = 0$ for $i \le C\log n$, and suppose that $|a_{C\log
    n + 1}| \ge \frac{1}{2}$. Recall that by Lemma \ref{lem:tau_2}, we
  have
  \[ \E\left[ \varphi(|g(\tau) - k_*|) \,\middle|\, F^k \right] \le O(C^{2/3} \log^{1/3} n). \]
  Then, Lemma \ref{lem:distinguishing-with-shifts} gives some index $j
  = O(C \log n)$ for which
  \begin{align*}
    \Big| \E\left[ \E_{\theta^S(\a)}(\wt{a}_j) \right] \Big| &\ge e^{- O(C^{1/3} \log^{1/3} n) - O(C^{2/3} \log^{1/3} n) } \\
    &= e^{-O(C^{2/3} \log^{1/3} n)}.
  \end{align*}
  Substituting into \eqref{eq:lp-approx1}, this is a contradiction of
  \eqref{eq:v-bound-hypothesis} for large enough $C$. We conclude that
  if \eqref{eq:v-bound-hypothesis} holds, then we must have $|w_0 -
  x_{k+1}| = |a_{C\log n + 1}| > \frac{1}{2}$.

  For the second claim, let $\w = \x(k+1:k+2C^2 \log n)$. Note that if
  on the event $F^k$ we have $\|V(\xt_{\w}) - V(\xt)\|_{C^2\log n,\infty}\neq 0$, then it means
  that there had to have been at least $C^2 \log n$ more deletions than
  insertions in the interval $[k, k + 2C^2 \log n]$. Thus,
  \[ \P( \{ \|V(\xt_{\w}) - V(\xt)\|_{C^2\log n,\infty}\neq 0\}  \cap F^k) \le e^{-\Omega(C^2 \log n)} = n^{-\Omega(C^2)}, \]
  and by Lemma \ref{lem:tau_2}, we have
  \[ \P( \|V(\xt_{\w}) - V(\xt)\|_{C^2\log n,\infty}\neq 0 \mid F^k) \le \frac{n^{-\Omega(C^2)}}{\P(F^k)} \le n^{-\Omega(C^2)}. \]
  Since the entries of $V(\xt)$ and $V(\xt_{\w})$ are all $0$ or $1$,
  this proves the second claim upon taking expectations.
\end{proof}

\subsection{Proof of Lemma \ref{lem:distinguishing-with-shifts}}

The remainder of this section is devoted to proving Lemma
\ref{lem:distinguishing-with-shifts}.

\begin{lemma} \label{lem:gen-func}
  Let $S$ be a bounded $\N$-valued random variable. Let
  $\a=(a_0,a_1,\dots)\in[-1,1]^{\N}$, and let $\at$ be the output from
  the insertion-deletion channel with deletion (resp.\ insertion)
  probability $q$ (resp.\ $q'$), applied to the randomly shifted
  string $\theta^{S}(\a)$. Let $\phi_1(w)=pw+q$, $\phi_2(w) =
  \frac{p'w}{1-q' w}$, and $\sigma(s)=\P[S=s]$ for $s\in\N$. Define
  \[ P(z) := \sum_{s=0}^{d} \sigma(s) z^s, \qquad Q(z):=\sum_{j= 0}^\infty a_j z^j. \]
  Then, for any $|w|<1$,
  \begin{equation} \label{eq:gen-func}
    \E\left[ \sum_{j\geq 0} \wt{a}_j w^j \right] =     
    \frac{pp'}{1-q'\phi_1(w)}
    \cdot
     P\left(\frac{1}{\phi_2\circ\phi_1(w)} \right ) \cdot Q(\phi_2\circ\phi_1(w)).
  \end{equation}
\end{lemma}
\begin{proof}
  Recall the construction of $\xt$ from $\x$ given in Section
  \ref{sec:reconstruction}, where we first insert a geometric number
  (minus one) bits before each bit of $\x$ and then delete each bit
  independently with probability $q$. From this description we see
  that we can sample $\wt\a$ by first setting
  $\wt\a^{(2)}=\theta^S(\a)$, then letting $\a^{(3)}$ be the string we
  get when sending $\a^{(2)}$ through the insertion channel with
  insertion probability $q'$ (and no deletions), and finally obtain
  $\wt\a$ by sending $\wt\a^{(3)}$ through the deletion channel with
  deletion probability $q$ (and no insertions).  Three elementary
  generating function manipulations (see, respectively, \cite[Lemma
    4.2]{PZ17}, \cite[Lemma 5.2]{NP16}, and \cite[Lemma 2.1]{NP16})
  give
  \[ \begin{split}
	&\E\left[ \sum_{j\geq 0} a_j^{(2)} w^j \right] = P(w^{-1})Q(w),\quad
	\E\left[ \sum_{j\geq 0} a_j^{(3)} w^j \,\,\bigg|\,\,\a^{(2)}\right] = \frac{\phi_2(w)}{w} \sum_{j\geq 0} a_j^{(2)} \phi_2(w)^j,\\
	&\E\left[ \sum_{j\geq 0} \wt a_j w^j\,\,\bigg|\,\,\a^{(3)} \right] = p\sum_{j\geq 0} a_j^{(3)} \phi_1(w)^j.
  \end{split} \]
  Combining these identifies we get \eqref{eq:gen-func}:
  \[ \begin{split}
	\E\left[ \sum_{j\geq 0} \wt a_j w^j \right]
	&= p\E\left[\sum_{j\geq 0} a_j^{(3)} \phi_1(w)^j\right]
	= p\frac{\phi_2\circ \phi_1(w)}{\phi_1(w)} \E\left[\sum_{j\geq 0} a_j^{(2)} \big(\phi_2\circ \phi_1(w)\big)^j\right]\\
	&=p\frac{\phi_2\circ \phi_1(w)}{\phi_1(w)}P\bigg(\frac{1}{\phi_2\circ \phi_1(w)}\bigg)Q(\phi_2\circ \phi_1(w)).
  \end{split} \]
\end{proof}

\def\cbe{C_{\mathrm{BE}}}

The following result is Corollary~3.2 of \cite{BE97} with $M=1$, $a =
\ell$ and $c_1 = \cbe$, observing that the class of polynomials whose
coefficients have modulus at most $1$ are in their class ${\mathcal
  K}_1^1$ and that their statement ${\mathcal K}_M := {\mathcal
  K}_M^0$ after their definition of ${\mathcal K}_M^\mu$ should be
ignored in favor of the correct statement ${\mathcal K}_M := {\mathcal
  K}_M^1$ occurring in their Corollary~3.2.

\begin{lemma}[Borwein and Erd{\'e}lyi 1997] \label{lem:BE}
  There is a universal constant $\cbe$ such that for any polynomial
  $f$ satisfying $|f(0)|=1$ and whose coefficients have modulus at
  most $1$, and for any arc $\alpha$ of the unit circle whose angular
  length is denoted $s \in (0,2\pi)$, we have
  \[ \sup_{z \in \alpha} |f(z)| \geq e^{-\cbe / s} \, .\]
\end{lemma}

\begin{proof}[Proof of Lemma \ref{lem:distinguishing-with-shifts}]
  Let
  \[ t^* = \argmin_{0 \le t \le n} \E\Big[ \varphi(|S - t|)\Big]. \]
  Define the quantities
  \begin{align*}
    \a &= \x^{(1)} - \x^{(2)} \in \{-1,0,1 \}^{\N}, \\
    L &= \max\left( n^{1/3}, \E\Big[ \varphi(|S - t^*|) \Big] \right), \\
    \rho &= 1 - 1/L^2.
  \end{align*}
  We first establish a general bound for M\"obius transformations
  appearing in Lemma \ref{lem:gen-func}.

  \begin{quote}
	\noindent{\underline{Claim:}} There is a constant $c_2 \in
    (0,1/20)$ depending only on $q,q'$ such that if $|\arg(z)| \leq
    c_2/L$, $|z|=1$, and $w = \phi_1^{-1} (\phi_2^{-1} (\rho \cdot
    z))$, then $|w|\leq 1-c_2/L^2$.

	\noindent{\underline{Proof:}} Observe that $\phi_2$
    (resp.\ $\phi_1$) is a M\"{o}bius transformation mapping $\D$ to a
    smaller disk which is contained in $\overline{\D}$, which is
    tangent to $\partial\D$ at 1, and which maps $\R$ to $\R$.  In
    particular, defining $\Psi := \phi_1^{-1}\circ\phi_2^{-1}$, we get
    by linearizing the map around $z=1$ that $\Psi(1+\wt z) = 1 + a
    \wt z + O(|\wt z|^2)$ for $a>1$ depending only on $q,q'$.  Writing
    $z=e^{i\theta}$, we have
	\begin{eqnarray*}
	  w & = & \Psi(\rho e^{i\theta})
	  \\ & = & 
	  1 + a(  \rho e^{i\theta} - 1) + O(|\rho e^{i\theta} - 1|^2)\\
	  & = & 1 + a\big( (1-L^{-2})(1+i\theta) - 1 \big) + O( \theta^2+L^{-4} ) \\
	  & = & 1 + a( -L^{-2} + i\theta ) + O( \theta^2+L^{-4} ),
	\end{eqnarray*}
	so $|w|< 1-c_2/L^2$ when $c_2=c_2 (q,q')$ is sufficiently small, and
	the claim is proved.
  \end{quote}

  Let $P$ and $Q$ be as in Lemma \ref{lem:gen-func}, and write
  $\wt{Q}(z) := z^{-n-1}Q(z)$, where our assumption that $\x$ and
  $\x'$ first differ in the $(n+1)$-th bit implies that $\wt{Q}(z)$ is
  a polynomial with $|\wt{Q}(0)| \geq 1/2$.

  Observe that $z \mapsto \wt{Q} (\rho \cdot z)$ has coefficients of
  modulus at most $1$, so we may apply Lemma~\ref{lem:BE} to find $z_0
  = e^{i\theta}$ with $|\theta| \leq c_2 / L$ such that $|\wt{Q} (\rho
  z_0)| \geq e^{-\cbe L/c_2}$. By definition of $c_2$, we see that
  $w_0 := \Psi (\rho \cdot z_0)$ satisfies $|w_0| \leq 1-c_2/L^2$. An
  illustration of the points $z_0$ and $w_0$ is given in
  Figure~\ref{fig:z_0}.

  \begin{figure}
	\centering
	\includegraphics[scale=1]{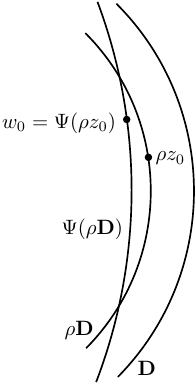} \\[2ex]
	\caption{Illustration of the points $z_0,w_0\in\C$ defined in the
      proof of Lemma~\ref{lem:distinguishing-with-shifts}.  We first choose $z_0 =
      e^{i\theta}$ for $|\theta|\leq c_2/L$, such that $|\wt Q(\rho
      \cdot z_0)|$ is bounded from below.  Then we observe that
      $|w_0|<1-c_2/L^2$, which helps us to bound the modulus of
      $\E[\sum_{j\geq 0} \wt a_j w_0^j]$ from below.}
	\label{fig:z_0}
  \end{figure}

  We next show that
  \begin{equation} \label{eq:P}
	\left | P \left ( \frac{1}{\rho z_0} \right ) \right | \geq \frac{1}{2} \, .
  \end{equation}
  To see this, define $\wt{P} (z) = z^{-t^*} P(z)$, which is an
  analytic function in the right half-plane. For all $z$ in the right
  half-plane satisfying $1 \le |z| \le \rho^{-1}$, differentiating
  $\wt{P}$ gives
  \begin{align*}
	|\wt{P}' (z)| &= \left| \sum_{j = 0}^{d} (j - t^*) \sigma(j) z^{j - t^* - 1} \right|
	\le \sum_{j = 0}^{d} |j - t^*|\cdot \sigma(j) \cdot |z|^{j - t^* - 1} \\
	&\le \E\Big[ |S - t^*| \cdot \rho^{-|S - t^*|} \Big] \le 4\E\Big[ |S - t^*| \cdot e^{|S - t^*|/L^2} \Big] \\
    &\le 4\E\Big[ |S - t^*| \cdot e^{|S - t^*|/n^{2/3}} \Big] = 4\E\Big[ \varphi(|S - t^*|) \Big] \le 4L.
  \end{align*}
  We also have
  \begin{align*}
	|\rho^{-1} z_0^{-1} - 1| &= \rho^{-1} |1 - \rho z_0| \le |z_0 - 1|
	+ \rho^{-1}(1 - \rho) \le \frac{c_2}{L}+ \frac{2}{L^2}.
  \end{align*}
  Therefore, for all sufficiently large $m$,
  \begin{align*}
	\left | P(\rho^{-1}z_0^{-1}) \right |
	&= \rho^{-\E S} \left | \wt{P}(\rho^{-1}z_0^{-1}) \right |
	\ge 1 - |\wt{P}(\rho^{-1}z_0^{-1}) - 1| \\
	&= 1 - \left| \int_1^{\rho^{-1}z_0^{-1}} \wt{P}'(z) \,dz \right|
	\ge 1 - |\rho^{-1}z_0^{-1} - 1| \cdot 4L \\
	&\ge 1-\left( \frac{c_2}{L}+ \frac{2}{L^2} \right)\cdot 4L \ge \frac{1}{2} \, ,
  \end{align*}
  proving~\eqref{eq:P}.

  Also note that the following quantity is bounded from below by a constant depending only on $p$ and $p'$:
  \begin{align*}
  \frac{pp'}{|1-q\phi_1(w_0)|}.
  \end{align*}	

  \def\csep{C_{\rm sep}}
  \def\cfwd{C_{\rm fwd}}
  Using Lemma \ref{lem:gen-func} and the above estimates, it follows that
  \begin{eqnarray} \label{eq:distinguish1}
	\left | \E \left [ \sum_{j\geq 0} \wt a_j w_0^j \right] \right |
	& \geq & \frac{pp'}{|1-q\phi_1(w_0)|}\cdot \left | P \left ( \frac{1}{\rho z_0} \right ) \right |
	\cdot  \rho^{n+1} \cdot|\wt{Q} (\rho z_0)|\\
	& \geq & \frac{pp'}{|1-q\phi_1(w_0)|}\cdot \frac{1}{2} \left ( 1 - \frac{1}{L^2} \right )^{n+1} e^{- \cbe L/c_2}
	\nonumber \\
	& \geq & e^{- O(L) - O(n/L^2)} = e^{-O(L)}. \nonumber
  \end{eqnarray}
  Since $|w_0|\leq 1-c_2/L^2$, for any $\cfwd > 1$,
  \begin{equation} \label{eq:distinguish2}
    \left|\sum_{j\geq \cfwd n} \E[\wt a_j] w_0^j\right|
    \leq
    \left | \sum_{j \geq \cfwd n} \left ( 1 - \frac{c_2}{L^2} \right )^j \right |
    \leq L^2 c_2^{-1} e^{- \cfwd L/c_2} = e^{-\Omega(\cfwd L)}.
  \end{equation}
  Combining \eqref{eq:distinguish1} and \eqref{eq:distinguish2}, by
  taking $\cfwd$ sufficiently large (depending only on $q$ and $q'$),
  we have
  \[ \begin{split}
	\E\left[ \sum_{j=0}^{\lceil \cfwd n \rceil-1} \left|\wt a_j w_0^j \right|\right]
	\geq
	\left|\E\left[ \sum_{j=0}^{\lceil \cfwd n \rceil-1} \wt a_j w_0^j \right]\right|
	\geq e^{-O(L)}.
  \end{split} \]
  It follows by the pigeonhole principle that there is a $j < \cfwd n$
  for which
  \[ |\E[\wt{a}_j]| \geq |\E[\wt{a}_j] w_0^j| \geq (2 \lceil\cfwd n\rceil)^{-1} e^{-O(L)} = e^{-O(L)}, \]
  as desired.
\end{proof}

\section{Proof of Theorem \ref{thm:main}}
\label{sec:proof-conclusion}
For each $k$, we will describe how to reconstruct $x_{k+1}$ assuming
we know $\x(0:k)$. This will involve applying the rules $\tau^k_1$ and
$\tau^k_2$. For purposes of analyzing the computational cost, we also
assume that the results of $\tau^i_1$ and $\tau^i_2$ have been saved
for all $i < k$. We assume throughout the proof that $k\geq 2C\log n$; the case $k<2C\log n$ is easier since we can reconstruct the first $2C\log n$ bits directly by the results of Section \ref{sec:bit-test} (with no shift).

We will use $n^{o(1)}$ traces to do the reconstruction of
$x_{k+1}$ with success probability at least $1 - n^{-2}$. Thus,
even if we reuse the sampled traces at each step, by a union bound,
reconstruction will succeed for all bits with high probability. The
computational cost for reconstructing this bit will also be $n^{o(1)}$
with probability at least $1 - n^{-2}$, so that the overall
computational cost is $n^{1 + o(1)}$ with high probability.

\subsection{Computing $\tau^k_1$ and $\tau^k_2$ in $n^{o(1)}$ time}

For each trace $\xt$, we want to first find $\tau^k_1$ and then find $\tau^k_2$. If we wanted to determine $\tau^k_1$ with probability 1 we would need to perform a sliding window of tests $T_{\ell,\lambda}(\x(k-\ell+1:k), \wt{\x}(k'-\ell+1:k'))$ where $k'$ potentially ranges from $\ell-1$ to $2n$. Each test only takes $\polylog$ time to perform, but as described, we are performing
$\Theta(n)$ tests, which exceeds our goal of $n^{o(1)}$.

To save on computation, we will only find estimates for $\tau^k_1$ and $\tau^k_2$, which are correct with high probability. Observe that if previously we had $\tau^i_1$
close to $f(i)$ for some $i$, then to find a (non-spurious) match, we
need only test bits in $\xt$ close to or after position $f(i)$. To carry out this
argument formally, we begin with the following definition.

\begin{definition}
  Suppose that $\x \not\in \bad$ and $\xt$ is a trace from $\x$. We
  say that $\xt$ is \textbf{progressively alignable} if the following
  hold:
  \begin{enumerate}
  \item For any $i \in \{1, \dots, n\}$ such that $\tau^i_1 < \infty$,
    we have $d(i, \tau^i_1) \le \log^2 n$ (recall \eqref{eq:d}).
  \item For any $t$ with $0 \le t \le n - e^{\log^{1/2} n}$, there is some $i
    \in [t, t + e^{\log^{1/2} n}]$ such that $\tau^i_1 < \infty$.
  \end{enumerate}
\end{definition}

The thresholds in the above definition have been chosen so that traces
are progressively alignable with very high probability, as shown by
the next lemma.

\begin{lemma} \label{lem:progressively-alignable}
  Suppose that $\x \not\in \bad$ and $\xt$ is a trace from $\x$. Then
  $\xt$ is progressively alignable with probability at least $1 -
  O(n^{-2})$.
\end{lemma}
\begin{proof}
  Let $\ell = C \log^{5/3} n$ and $\lambda = C \log^{2/3} n$ as in
  Lemma \ref{lem:tau_1}. To see that the first property occurs with
  probability $1 - O(n^{-2})$, we may use the same argument as in the
  proof of Lemma \ref{lem:tau_1}, with the only modification being
  that in the definitions of the events $E'$ and $F$, the quantities
  $\frac{1}{10} C \log n$ and $\frac{1}{20} C \log n$ should be
  changed to $\log^2 n$ and $\frac{1}{2} \log^2 n$, respectively.

  For the second property, note that since $\x$ is coarsely
  well-behaved, each of its substrings $\w$ of length $\ell$ exhibits
  robust bias at scale $\lambda$. Thus, we may apply Lemma
  \ref{lem:T-match-lower} to conclude that with probability at least
  $e^{-O(\ell/\lambda^2)} \ge e^{-O(\log^{1/3} n)}$, the part of
  the trace coming from $\w$ contains a match for the test
  $T_{\ell,\lambda}(\w, \,\cdot\,)$.

  In any interval $[t, t + e^{\log^{1/2} n}]$, we have at least
  $e^{\Omega(\log^{1/2} n)}$ disjoint intervals of length $\ell$, and
  each of these has independently a $e^{-O(\log^{1/3} n)}$ chance of
  producing a match. The chance of not having a single match over this
  whole interval is therefore at most
  \[ \left( 1 - e^{-O(\log^{1/3} n)} \right)^{e^{\Omega(\log^{1/2} n)}} = \exp\left( -\frac{e^{\Omega(\log^{1/2} n)}}{e^{O(\log^{1/3} n)}} \right) \le n^{-3}. \]
  Taking a union bound over all $t$ establishes the second property.
\end{proof}

Let us now specify our algorithm for computing (an estimate of)
$\tau^k_1$, which is only guaranteed to give the right value (i.e.,
the value defined in Lemma \ref{lem:tau_1}) if the trace is
progressively alignable, but also costs only $n^{o(1)}$
operations. The algorithm is to first look for $i \in [ k -
  2e^{\log^{1/2} n}, k - e^{\log^{1/2} n} ]$ such that $\tau^i_1
< \infty$. Then, we evaluate $\tau^k_1$ as
\[ \inf \left\{ \tau^i_1 \le k' \le \tau^i_1 + 3e^{\log^{1/2} n} : T_{\ell,\lambda}(\x(k-\ell:k), \xt(k'-\ell:k')) = 1 \right\}. \]
As long as the trace is progressively alignable and our stored values
of $\tau^i_1$ are correct, this gives us the right answer, i.e., our
estimate for $\tau_1^k$ is equal to the true value of
$\tau_1^k$. The above test takes only $e^{O(\log^{1/2} n)} =
n^{o(1)}$ operations, and by Lemma \ref{lem:progressively-alignable}
and a union bound, the overall probability of getting at least one wrong result is less than $n^{-1}$.

Next, in order to calculate $\tau^k_2$, it is necessary to identify
the ``good'' interval $I$ in Definition \ref{def:finely-behaved}. For
a given interval $I$, it can be easily checked whether $\x(I)$ has
robust bias at scale $\lambda$, but it is not as straightforward to
explicitly calculate
\[ \P_{\x}\left( \cQ_{\ell,\lambda}(I, J) \right). \]
However, we can estimate this probability to high accuracy by
Monte-Carlo simulation. Since the relevant interval $J$ is only of
logarithmic size, each simulated sample can be produced in $n^{o(1)}$
time. By Hoeffding's inequality, $e^{O(C^2 \log^{1/3} n)}$ samples are
enough get an accuracy of $e^{-\Omega(C^2 \log^{1/3} n)}$ with
probability $1 - n^{-2}$. This level of accuracy is small compared to
the probability bound in Definition \ref{def:finely-behaved}, so it is
accurate enough for all of our analysis to carry through.

\subsection{Determining the next bit}

Our algorithm will work by sampling $N := e^{C^2 \log^{1/3} n}$
traces. Let $N_1$ denote the number of traces for which $\tau^k_2 <
\infty$, and let us number these traces $\xt^{(1)}, \ldots ,
\xt^{(N_1)}$. Note that by Lemma \ref{lem:tau_2}, we have $N_1 =
e^{\Omega(C^2 \log^{1/3} n)}$ with probability at least $1 -
e^{-\Omega\left(e^{\Omega(\log^{1/3} n)}\right)} \ge 1 - n^{-2}$.

Our reconstruction strategy is based on Lemma \ref{lem:approx-lp}. Let
$m := C^2 \log n$. Recall from \eqref{eq:maxm} the notation $\Linfm{\a}$ for any $\a \in [0,1]^\N$. 
Suppose we are able to (approximately) solve the following
minimization problem (where $v$ is as in Definition \ref{def:v}):
\begin{equation} \label{eq:lp}
  \min_{\w \in [0,1]^{2m}} \Linfm{v(\w) - v(\x(k+1:\infty))}.
\end{equation}
Then, by Lemma \ref{lem:approx-lp}, we could recover $x_{k+1}$ by
rounding $w_0$ to the nearest integer (either $0$ or $1$). Note that
\eqref{eq:lp} is a linear program with $O(\log n)$ variables and
constraints, and so it can be solved in $\polylog$ time by
e.g.\ \cite{karmarkar}.\footnote{More precisely, in \cite[Section
    1.6]{karmarkar} it is proved that if $L$ is the number of bits in
  the input and $2m$ is the number of variables then the problem can
  be solved in time $O(m^{3.5}L^2\log L\log\log L)$. We can round the
  coefficients of $v$ to the nearest multiple of
  $m^{-1}e^{-K\log^{1/3}n}=e^{-\Theta(K\log^{1/3}n)}$. For $K\gg 1$,
  this gives $L=O( m K\log^{1/3}n)$, so we can find a solution
  which differs from the optimal solution by
  $e^{-\Omega(K\log^{1/3}n)}$ and has running time of order
  $m^{3.5}L^2\log L\log\log L \ll \log^7 n$.}

However, two issues arise: we do not have
direct access to the quantity $v(\x(k+1:\infty))$, nor are we able to
evaluate quantities like $v(\w)$ directly.

To address the first issue, we can estimate $v(\x(k+1:\infty))$ using our
traces. Consider the empirical mean
\[ \wh{v} = \frac{1}{N_1} \sum_{i = 1}^{N_1} V(\xt^{(i)}). \]
This is a sum of i.i.d.\ vectors with entries in $[0, 1]$ whose
expectation is the desired vector $v(\x(k+1:\infty))$. Thus, by
Hoeffding's inequality, we have with probability at least $1 - n^{-2}$
that
\[ \Linfm{\wh{v} - v(\x(k+1:\infty))} \le e^{-\Omega(C^2 \log^{1/3} n)}. \]

To address the second issue, we will estimate $v(\w)$ by Monte-Carlo
simulation. Let $\e_i$ denote the vector with $1$ in the $i$-th entry
and $0$ elsewhere; we first estimate the quantities $v(\e_i)$. Since
we already know $\x(0:k)$, we can simulate drawing a trace from
$\x_{\e_i}$ (Definition \ref{def:v}) and compute $\tau^k_2$. However, this once again takes
$O(n)$ time, because we have to scan through the whole string.

Instead, we sample a trace $\xt'_{\e_i}$ from the shortened string
$\x_{\e_i}(k - \log^2 n :\infty)$ and evaluate $V(\xt'_{\e_i})$. Note that
the trace $\xt'_{\e_i}$ is equivalent to removing the first $f(k -
\log^2 n)$ bits from a trace $\xt_{\e_i}$ of the full string
$\x_{\e_i}$. Coupling $\xt'_{\e_i}$ and $\xt_{\e_i}$ in this way,
$V(\xt'_{\e_i})$ is usually the exact same as $V(\xt_{\e_i})$; the
only way they can differ is if $g(\tau^k_1(\xt_{\e_i})-\ell) \le k - \log^2
n$, which by Lemma \ref{lem:tau_1} happens with probability
$O(n^{-2})$. Note that we do not know that $\x_{\e_i} \not\in \bad$ (in fact, we typically have $\x_{\e_i}\in\bad$). However, we have $\x_{\e_i}(0:k)=\x(0:k)$, and the initial part of the string is the most relevant part when we do the alignment, since the string we use to align is chosen as a substring of $\x_{\e_i}(0:k)=\x(0:k)$. Furthermore, adding the string $\e_i$ at the end will cause false positives with very small probability since the bit statistics of this string are very different from those of the string we use to align. 
It follows that
\begin{align*}
&\Linfm{\E[ V(\xt'_{\e_i}) \mid \tau^k_2(\xt'_{\e_i}) < \infty ] - \E[ V(\xt_{\e_i}) \mid \tau^k_2(\xt_{\e_i}) < \infty ]} \\
&\qquad = \Linfm{ 
	\frac{\E[ V(\xt_{\e_i}) \1_{\tau^k_2(\xt_{\e_i}) < \infty} ]+O(n^{-2})}{\P[\tau^k_2(\xt_{\e_i}) < \infty]+O(n^{-2}) } 
	- 
	\frac{\E[ V(\xt_{\e_i}) \1_{\tau^k_2(\xt_{\e_i}) < \infty} ]}{\P[\tau^k_2(\xt_{\e_i}) < \infty]}}
= O(n^{-1.5}).
\end{align*} 
We can also see from Lemma \ref{lem:tau_2} that
\begin{align*}
  &\Linfm{v(\e_i) - \E\left[ V(\xt_{\e_i}) \,\middle|\, \tau^k_2(\xt_{\e_i}) < \infty \right]} \\
  &\qquad= \Linfm{\E\left[ V(\xt_{\e_i}) \,\middle|\, F^k \right] - \E\left[ V(\xt_{\e_i}) \,\middle|\, \tau^k_2(\xt_{\e_i}) < \infty \right]} \le O(n^{-1.5}).
\end{align*}
Thus, by performing this simulation $e^{C^2 \log^{1/3} n}$ times and
averaging the results, we are able to provide an estimate $v'(\e_i)$
of $v(\e_i)$ to within $e^{-\Omega(C^2 \log^{1/3} n)}$ error with
probability at least $1 - O(n^{-2})$. We can extend this linearly to
all possible inputs $\w \in [0, 1]^{2m}$ by setting
\[ v'(\w) = \sum_{i = 1}^{2m} w_i v'(\e_i), \]
and we see that overall $\Linfm{v'(\w) - v(\w)} \le e^{-\Omega(C^2
  \log^{1/3} n)}$.

We can then solve the modified optimization problem
\[ \min_{\w \in [0,1]^{2m}} \Linfm{v'(\w) - \wh{v}}, \]
which is an approximation of our original problem. With probability at
least $1 - n^{-2}$, the objective function in the above problem is
within $e^{-\Omega(C^2 \log^{1/3} n)}$ of the objective in
\eqref{eq:lp}. In this case, as long as $C$ is large enough, our
minimizer $\w^*$ will satisfy the hypothesis of Lemma
\ref{lem:approx-lp}, and so we can correctly extract the next bit
$x_{k+1}$ as the closer of $0$ or $1$ to $w^*_1$. This completes our
analysis and establishes Theorem \ref{thm:main}.

\vskip10pt 

{\bf Acknowledgements} We thank Margalit Glasgow and the anonymous referees for helpful comments.

\bibliographystyle{hmralphaabbrv}
\bibliography{ref}

\end{document}